\documentclass[final,leqno]{siamltex}

\usepackage{rotating}
 \usepackage{amsmath}
\usepackage{graphics}
 \usepackage{amsfonts}
 \usepackage{amssymb}
 \usepackage{amscd}
\usepackage{colordvi}
\usepackage{color}
\usepackage{pdfsync}
\usepackage{url}
\usepackage{enumerate}
\usepackage{booktabs}
\usepackage{multicol}
\usepackage[english]{babel}
\usepackage{blindtext}
\usepackage[normalem]{ulem} 
\usepackage{subfigure} 

\newtheorem{thm}{Theorem}[section]
\newtheorem{cor}[thm]{Corollary}
\newtheorem{lem}[thm]{Lemma}
\newtheorem{prop}[thm]{Proposition}
\newtheorem{ex}[thm]{Example}

\newtheorem{rem}[thm]{Remark}

\def \R { {\mathbb R} }

\def \Z { {\mathbb Z} }

\title{The characteristic imset polytope of Bayesian networks with
ordered nodes}

\author{Jing Xi\thanks{Mathematics Department, North Carolina State
    University, 3600 Univeristy,  2108 SAS Hall,
2311 Stinson Drive
Raleigh, NC 27695-8205, U.S.A.
({\tt tykiallen@gmail.com}).}
        \and Ruriko Yoshida\thanks{Statistics Department, University of Kentucky, 
        325 Multidisplinary Science Building, 725 Rose Street
        Lexington, KY 40536-0082, U.S.A. ({\tt ruriko.yoshida@uky.edu})}}

\begin{document}

\maketitle

\begin{abstract}
In 2010, M. Studen\'y, R. Hemmecke, and S. Linder explored a new
algebraic description of graphical models, called characteristic
imsets. Compare with standard imsets, characteristic imsets have
several advantages: they are still unique vector representative of
conditional independence
structures, they are $0$-$1$ vectors, and they are more intuitive in terms
of graphs than standard imsets. After defining a characteristic imset
polytope (cim-polytope) as the convex hull of all characteristic imsets with a given
set of nodes, they also
showed that a model selection in graphical models, which maximizes a quality criterion, can be converted into a linear programming problem over the cim-polytope. However, in general, for a fixed set of nodes, the cim-polytope can have exponentially many vertices over an exponentially high dimension. 
Therefore, in this paper, we focus on the family of directed acyclic graphs
(DAGs) whose nodes have a fixed order.  This family includes
diagnosis models which can be
described by Bipartite graphs with a set of $m$ nodes and a set of $n$
nodes for any $m, n \in \Z_+$.    In this paper, 
we first consider
cim-polytopes for all diagnosis models and show that these polytopes
are direct products of simplices. Then we give a combinatorial description of all
edges and all facets of these polytopes. Finally, we generalize
these results to the cim-polytopes for all Bayesian networks with a
fixed underlying ordering of nodes with or without fixed (or forbidden) edges.  
\end{abstract}

\begin{keywords} 
 graphical model, characteristic imset
polytope, diagnosis model, Bipartite graph, directed acyclic graphs.
\end{keywords}

\begin{AMS}
51M20, 52B20, 13P25.
\end{AMS}

\pagestyle{myheadings}
\thispagestyle{plain}
\markboth{JING XI AND RURIKO YOSHIDA}{CIM POLYTOPES}

\section{Introduction}

Bayesian networks (BNs), also known as belief networks, Bayes
networks, Bayes(ian) models or probabilistic directed acyclic 
graphical models, find their applications to model knowledge in many
areas, such as computational biology and bioinformatics (gene regulatory networks,
protein structure, gene expression analysis \cite{expression} learning epistasis
from GWAS data sets \cite{epi}) and medicine \cite{meds}.
BNs are a part of the family 
of probabilistic graphical models (GMs). These
graphical structures represent knowledge
about probabilistic structures for a statistical model. More precisely, each node
in the graph represents a random variable and 
an edge between the nodes represents probabilistic
dependencies among the random variables corresponding to the nodes
adjacent to the edge \cite{Lauritzen}.
BNs correspond to GM structure known
as a directed acyclic graph (DAG) defined by the
set of nodes (vertices) and the set of directed edges.

In order to infer parameters from the observed data set,  we first apply
a model selection criterion called {\em quality criterion}, which
provides a way to construct highly 
predictive BN models from data by choosing the graph which gives the
given criteria, such as Bayesian Information Criteria (BIC) \cite{BIC} or Akaike
Information Criteria (AIC) \cite{akaike}, maximum (see \cite{studeny2008} for
more details on quality criterions).  Intuitively a
quality criterion is a function, $\mathcal{Q}(G, D)$, which takes a DAG,
$G$, and an observed data set, $D$, to evaluate how good the DAG $G$ to
explain the observed data $D$.  Note that different DAGS, $G_1, G_2$
may have the same conditional independences (CIs).  In that case we say
$G_1, G_2$ are {\em Markov equivalent}.  When researchers wish to infer
the CIs of the BN structure from the observed data set one represents
each set of Markov equivalent graphs by one graph called the {\em
  essential graph} the corresponding Markov equivalence class of DAGs \cite{Andersson}.
In this paper we focus on quality criterions $\mathcal{Q}(G, D)$, such
that $\mathcal{Q}(G_1, D) = \mathcal{Q}(G_2, D)$ if and only if $G_1,
G_2$ are Markov equivalent. 

Since in general there are super exponentially
many essential graphs with a fixed set of nodes $N$, maximizing the quality
criterion,  $\mathcal{Q}(G, D)$, over all possible essential graphs with $N$ is
known to be NP-hard.  Studen\'y developed an algebraic representation
of each essential graph $G$ called a {\em standard imset}, of $G$,
which is an integral vector representation of $G$ in $\R^{2^{|N|} - |N| -
  1}$.  From the view of this setting a criterion function
$\mathcal{Q}(G, D)$ is a dot product of vectors in  $\R^{2^{|N|} - |N| -
  1}$.  In 2010, M. Studen\'y, J. Vomlel, and 
R. Hemmecke showed that maximizing the $\mathcal{Q}(G, D)$ over all
essential graphs can be formulated as a linear programming problem
over the convex hull of standard imsets for all possible
essential graphs \cite{studeny2010}.  This gives us a systematic way
to find the best criterion with the optimality certificate rather than finding the best
criterion by the brute-force search.  Then 
M. Studen\'y, R. Hemmecke, and S. Linder 
explored an alternative
vector representative of the BN structure, called {\em characteristic
imsets}.  Compare with standard imsets, characteristic imsets have
several advantages: they are still unique vector representative of
conditional independence
structures; they are $0$-$1$ vectors; and they are more intuitive in terms
of graphs than standard imsets \cite{studeny2010a}. 

In general, however, the dimension of the convex hull of the
characteristic imsets with the fixed set of nodes $N$, called a {\em
  characteristic imset polytope} (cim-polytope), is
exponentially large and there are double exponentially many vertices (cim-polytope)
as well as 
facets of the cim-polytope.  Thus it is infeasible to
optimize by software if $|N| > 6$.   In order to solve the LP problem for a larger
$|N|$, we need to understand 
the structure of the cim-polytope, such as
combinatorial description of edges and
facets of
the polytope so that we might be able to apply a simplex method to
find an optimal solution.  However, in general, it is challenging because there
are too many facets and too many edges of the polytope.  Therefore
here we start with a particular family of BN models, namely {\em
  diagnosis models}.  

In medical studies, researchers are often interested in probabilistic
models in order for them to correctly diagnose a disease from a
patient symptoms.  The diagnoses models, also known as the Quick
Medical Reference (QMR) diagnostic model, is
introduced in \cite{Shwe} to diagnose a disease from a given set of
symptoms of a patient.
Therefore, here we focus on diagnosis models (e.g.,
\cite{Locas2001}). Under this model, a DAG representing the model is a
bipartite graph with two sets of nodes, one representing $m$ diseases and
one representing $n$ symptoms, and set of directed edges from nodes
representing diseases to nodes representing symptoms (see Figure
\ref{bipartite0} for an example).  

 In this paper, first, 
we are able to find an explicit combinatorial description of all
edges of the cim-polytopes for diagnosis models with
fixed $m$ and $n$, that is, if $G_1, \, G_2$ are graphs representing
two diagnosis models such that all symptoms have the same parents
in $G_1$ and in $G_2$ except one symptom, then the
characteristic imsets representing $G_1, \, G_2$ form an edge of the
cim-polytope for diagnosis models. Then
we prove that these cim-polytopes are direct products
of $n$ many $(2^m-1)$ dimensional simplices, and an explicit description of all facets of them can be given based on this structure.   Finally we generalize
 these results for the cim-polytopes for
BNs with a fixed underlying ordering with or without
fixed (or forbidden) edges. 

This paper is organized as follows. In Section \ref{prop} we introduce
notation, and we state some definitions as well as propositions and
their proofs.   Section \ref{BN:diagnosis} shows the description of the
cim-polytopes  for  diagnosis models and Section
\ref{BN:generalize} shows the description of the
cim-polytopes  for Bayesian networks with a fixed
underlying ordering. Proofs of some of properties, lemmas, and
theorems can be found at Section \ref{proof:diagnosis}, Section
\ref{proof:simple}, and Section \ref{proof:generalize}.  We end with a 
discussion of our future work in Section \ref{dis}.

\section{Definitions and propositions for diagnosis models}\label{prop}

\begin{figure}[!htp]
\begin{center}
\scalebox{0.8}{
\includegraphics{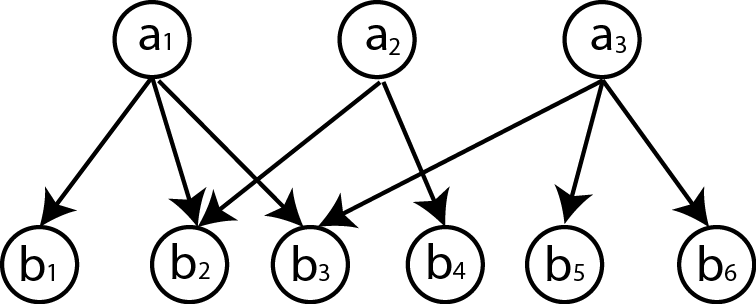}
}
\end{center}
\caption{An example of Bipartite Graph, $m=3$, $n=6$.}
\label{bipartite0}
\end{figure}

In this section we state some notation and remind readers some
definitions.  
\begin{definition}\label{def_bipar}
A Diagnosis Model can be described by a Bipartite Graph whose nodes $N=\{a_1,\ldots,a_m\} \cup \{b_1,\ldots,b_n\}$ can be divided into disjoint sets $A=\{a_1,\ldots,a_m\}$ and $B=\{b_1,\ldots,b_n\}$. Nodes in $A$ can be interpreted as diseases and nodes in $B$ can be interpreted as symptoms. Every single edge can only be drawn from a disease to a symptom. An example is given by Figure \ref{bipartite0}.\\
For fixed $A$ and $B$, where $|A|=m$ and $|B|=n$, we define notation: $\mathcal{G}_{m,n}=\{$All possible directed bipartite graphs defined in Definition \ref{def_bipar} based on $A$ and $B \}$. 
\end{definition}

Recall that we have the definition of Characteristic Imset.\\

\begin{definition}
Let $G$ be an acyclic directed graph over $N$. The {\em characteristic imset} for $G$ can be introduced as a zero-one vector $c_G$ with components $c_G(S)$ where $S\subseteq N$, $|S|\geq 2$ given by
\begin{center}
$c_G(S)=1\Longleftrightarrow \exists \ i\in S$ such that $j\in pa_G(i)$ for $\forall\ j\in S\backslash \{ i\}$
\end{center}
where $j\in pa_G(i)$ means $G$ includes the edge from $j $ to $i$.
\end{definition}

\begin{prop}\label{prop:diag1}
Fix $A=\{a_1,\ldots,a_m\}$ and $B=\{b_1,\ldots,b_n\}$. Assume $G \in \mathcal{G}_{m,n}$ and $|N| = m+n >2$. Then $c_G(T)$ is possible to take value $1$ if and only if $T$ has the form of $a_{i_1} \ldots a_{i_k} b_j$, where $1\leq k \leq m$, $\{i_1,\ldots,i_k\}\subseteq \{1,\ldots,m\}$ and
$j\in \{1,\ldots,n\}$.
\end{prop}

\begin{proof}
Notice that $\forall T\subseteq N$, $|T|\geq 2$, we can write T in the form of:
\begin{equation}
\begin{array}{rl}
T=a_{i_1} \ldots a_{i_k} b_{j_1} \ldots b_{j_l}, \mbox{ where}
 & 0\leq k \leq m, \{i_1,\ldots,i_k\}\subseteq \{1,\ldots,m\}, \\
 & 0\leq l \leq n, \{j_1,\ldots,j_l\}\subseteq \{1,\ldots,n\}, \\
 & k+l \geq 2.
\end{array}
\end{equation}
We need to prove that $l$ can neither be 0 nor greater than 1, i.e. $l=1$.
\begin{itemize}
\item[(a)] 
If $l=0$. $\forall\ s,t \in \{i_1,\ldots,i_k\}$, by Definition \ref{def_bipar}, $a_s \rightarrow a_t$ is not in $G$. This means $a_s \notin pa_G(a_t)$. Hence 
$\forall\ t \in \{i_1,\ldots,i_k\}$, $T\setminus \{a_t\} \nsubseteq pa_G(a_t)$. $c_G(T)=0$.
\item[(b)] 
If $l>1$. Similarly with above, by Definition \ref{def_bipar}, $\forall\ s',t' \in \{j_1,\ldots,j_l\}$, $b_{s'} \notin pa_G(b_{t'})$. Moreover, $\forall\ t \in \{i_1,\ldots,i_k\}$
and $t' \in \{j_1,\ldots,j_l\}$, $b_{t'} \notin pa_G(a_t)$. $c_G(T)=0$.
\end{itemize}
\end{proof}

\begin{prop}\label{prop:diag2}
Notation is adopted from Proposition \ref{prop:diag1}. Suppose $T$ has the form of  $a_{i_1} \ldots a_{i_k} b_j$, where $1\leq k \leq m$, $\{i_1,\ldots,i_k\}\subseteq \{1,\ldots,m\}$ and $j\in \{1,\ldots,n\}$, then $c_G(T)=\prod_{s=i_1,\ldots,i_k} c_G(a_s b_j)$.
\end{prop}

\begin{proof}
Again by Definition \ref{def_bipar}, $\forall\ s,t \in \{i_1,\ldots,i_k\}$, $a_s \notin pa_G(a_t)$. Therefore:
\begin{equation}
\begin{array}{rl}
c_G(T)=1 & \Longleftrightarrow \{a_{i_1} \ldots a_{i_k} \} \subseteq pa_G(b_j) \\
& \Longleftrightarrow a_s\in pa_G(b_j), \ \forall s=i_1,\ldots,i_k\\
& \Longleftrightarrow c_G(a_s b_j)=1, \ \forall s=i_1,\ldots,i_k.
\end{array}
\end{equation}
Recall that $c_G(T)$ is binary. Thus $c_G(T)=\prod_{s=i_1,\ldots,i_k} c_G(a_s b_j)$.
\end{proof}

\begin{rem}\label{rem:diag1}
Proposition \ref{prop:diag2} implies that $\forall G \in \mathcal{G}_{m,n}$, $c_G$ is determined by only $m\cdot n$ coordinates, $\{c_G(a_i b_j): i = 1,\ldots,m, j = 1,\ldots,n\}$, i.e. the existence of directed edges $a_i \rightarrow b_j$, $i = 1,\ldots,m$ and $j = 1,\ldots,n$. Another way to see this property is that $\forall\ G \in \mathcal{G}_{m,n}$, $G$ can be determined by $pa_G(b_j)$, $b_j\in B$. Thus if we consider a permutation of coordinates in $c_G$ that corresponds to a permutation of $T$ where $T$ has the form in Proposition \ref{prop:diag1}, then these coordinates can be broken into $n$ parts:
\[
\underline{a_1 b_1,\ldots, a_m b_1,\ldots, a_1 \ldots a_m b_1}, 
\uwave{a_1 b_2,\ldots, a_m b_2, \ldots, a_1 \ldots a_m b_2}, \ldots, \underline{a_1 b_n,\ldots, a_1 \ldots a_m b_n},
\]
where the $s$-th part of coordinations $c_G(T)$, $T\in \{a_1 b_s,\ldots, a_m b_s, a_1 a_2 b_s, \ldots, a_1 \ldots a_m b_s\}$ only depend on $pa_G(b_s)$, and different parts are completely irrelevant in the sense that $pa_G(b_s)$, $b_s \in B$, can be decided separately.
\end{rem}

\begin{prop}\label{prop:diag3}
Fix $m$ and $n$. The number of elements in $\mathcal{G}_{m,n}$ is $2^{mn}$.
\end{prop}

\begin{proof}
This is trivial because of Remark \ref{rem:diag1} since there are $mn$ possible edges that can be assigned: $a_i \rightarrow b_j$, where $i = 1,\ldots,m$ and $j = 1,\ldots,n$, and there are  
$\sum_{k = 0}^{mn} {mn \choose k} = 2^{mn}$
many possible ways to assign the existence of these edges. 
\end{proof}

\begin{prop}\label{prop:diag4}
Suppose $G \in \mathcal{G}_{m,n}$. The number of non-zero coordinates in $c_G$ is at most $n\cdot (2^m-1)$.
\end{prop}

\begin{proof}
This result is straightforward from Proposition \ref{prop:diag1} by counting the number of coordinates $c_G(T)$, where $T$ has the form shown in Proposition \ref{prop:diag1}. Note that when $|T|> m+1$, $\exists \ b_{j_1},b_{j_2}\in \{1,\ldots,n\}\ s.t. \ b_{j_1},b_{j_2}\in T$, i.e. $c_G(T) = 0$ by Proposition \ref{prop:diag1}. When $2 \leq |T| \leq m+1$, the number of coordinates of form $c_G(a_{i_1}\ldots a_{i_{|T|-1}} b_j)$, where $\{i_1,\ldots,i_{|T|-1}\}\subseteq \{1,\ldots,m\}$ and $j\in \{1,\ldots,n\}$, is
${m \choose |T|-1} \cdot n$.
Hence the number of possible non-zero coordinates is:
\[
\sum\limits_{|T|=2}^{m+1} {m \choose |T|-1} \cdot n
=n\cdot \sum\limits_{k=1}^{m} {m \choose k} 
=n\cdot (2^m-1).
\]
\end{proof}

\begin{definition}
Recall several definitions in elementary geometry (see \cite{Ziegler} for more details on polyhedral geometry):
\begin{itemize}
\item
a {\em closed convex polyhedron} (which will be indicated as {\em polyhedron} for short) in $\mathbb{R}^q$ can be defined by a system of linear inequalities:
\begin{center}
$\{\mathbf{x}\in \mathbb{R}^q : A\mathbf{x} \leq \mathbf{b} \}$
\end{center}
where $A$ is a $p\times q$ matrix in $\mathbb{R}^{p\times q}$ and $\mathbf{b}$ is a vector in $\mathbb{R}^{p}$;
\item
a {\em closed convex polytope} (which will be indicated as {\em polytope} for short) is defined as the convex hull of a finite set of points;
\item
if a polyhedron is bounded, then it is a polytope;
\item
for a polytope $\mathbf P$, we define $vert(\mathbf P)$ as the set of vertices of $\mathbf P$;
\item
A {\em d-simplex} is a d-dimensional polytope which has exactly $d+1$
vertices. It is notated as $\Delta_d$.
\end{itemize}
Let $DAGs(N)$ be the set of all directed acyclic graphs over $N$, and consider a class of graphs $\mathcal G \subseteq DAGs(N)$ that contains all graphs which we are interested in.
We call the convex hull of $\{c_G : G \in \mathcal G\}$, $\mathbf
P_{\mathcal G} = conv\{c_G : G \in \mathcal G\}$ the {\em
  characteristic imset polytope} (cim-polytope) for $\mathcal G$. Note that it is obvious that $vert(\mathbf P_{\mathcal G}) = \{c_G : G \in \mathcal G\}$.
\end{definition}

For fixed $A$ and $B$ in Definition \ref{def_bipar}, define $\mathbf P_{m, n} := \mathbf P_{\mathcal{G}_{m,n}}$. Proposition \ref{prop:diag4} implies that the dimension of $\mathbf P_{m, n}$ is at most $n\cdot (2^m-1)$. We will show that the dimension of $\mathbf P_{m, n}$ is actually exactly $n\cdot (2^m-1)$.

\section{The cim-polytopes for diagnosis models}\label{BN:diagnosis}

\subsection{Combinatorial description of edges on $\mathbf{P}_{m,n}$}\label{BN:diagedges}

\begin{definition}
Consider a class of graphs $\mathcal G \subseteq DAGs(N)$. $\forall
G$, $H\in \mathcal{G}$, $G$ and $H$ are called {\em neighbors} if $c_G$ and $c_H$ form an edge in $\mathbf P_{\mathcal G}$, the cim-polytope for $\mathcal G$.
\end{definition}

\begin{lem}\label{lem:onesymp}
Fix $m$. Suppose $G_1$, $G_2 \in\mathcal{G}_{m,1}$ are arbitrary two distinct graphs in $\mathcal{G}_{m,1}$. Then $G_1$ and $G_2$ are neighbors, i.e. $c_{G_1}$ and $c_{G_2}$ form an edge in $\mathbf P_{m,1}$.
\end{lem}

\begin{proof}
See Section \ref{proof:diagnosis}.
\end{proof}

\begin{thm}\label{thm:diagedge}
Fix $m$ and $n$. Two graphs, $G_1$, $G_2\in \mathcal{G}_{m,n}$ are neighbors if and only if $\exists\ b_i\in B$ such that $pa_{G_1}(b_i)\neq pa_{G_2}(b_i)$ and $pa_{G_1}(b_j) = pa_{G_2}(b_j)$, $\forall\ b_j\in B$ and $b_j\neq b_i$, i.e. all nodes but one have exactly the same parent sets in $G_1$ and $G_2$.
\end{thm}

\begin{proof}
See Section \ref{proof:diagnosis}.
\end{proof}

\subsection{$\mathbf{P}_{m,n}$ is a direct product of simplices}\label{BN:diagfacets}

\begin{thm}\label{thm:diagnedge}
Fix $m$ and $n$. For an arbitrary $G\in \mathcal{G}_{m,n}$, $G$ has $n\cdot (2^m-1)$ many neighbors.
\end{thm}

\begin{proof} 
See Section \ref{proof:diagnosis}.
\end{proof}

\begin{rem}\label{rem:diag2}
Theorem \ref{thm:diagnedge} implies that every vertex of $\mathbf P_{m,1}$ has $(2^m-1)$ neighbors. Since $|vert(\mathbf P_{m,1})|=2^m$ (by Proposition \ref{prop:diag3}), $\mathbf P_{m,1}$ is a simplex with dimension $(2^m-1)$, i.e. $\mathbf P_{m,1} = \Delta_{2^m-1}$.
\end{rem}

\begin{thm}\label{thm:diagprod}
$\mathbf{P}_{m,n}$ is the direct product of $n$ many $\Delta_{2^m-1}$, i.e.
\begin{center}
$\mathbf{P}_{m,n}= \underbrace{\Delta_{2^m-1} \times \Delta_{2^m-1} \times \cdots \times \Delta_{2^m-1}}_{\mbox{n many}}$.
\end{center}
And the $i_{th}$ simplex is $\mathbf{P}_{m,1}$ with the same diseases $A$ and only one symptom $\{b_i\}$.
\end{thm}

\begin{proof} Fix $m$, we are going to prove the equality by induction on $n$.
\begin{itemize}
\item
$n=1$. See Remark \ref{rem:diag2};
\item
Fix $q\in \mathbb{Z}^+$. Suppose the equality holds for $\mathbf{P}_{m,n}$, $\forall \ n<q$, then we need to prove that it also holds for $\mathbf{P}_{m,q}$. Recall that for $\mathcal{G}_{m,q}$, the symptoms are: $B=\{b_1,b_2,\ldots,b_q\}$.

First, we need to prove: $\mathbf{P}_{m,q}\subseteq \mathbf{P}_{m,q-1}\times \mathbf{P}_{m,1}$.

Similarly with the proof of Theorem \ref{thm:diagedge}, $\forall\ G\in\mathcal{G}_{m,q}$, we define graphs:
\begin{itemize}
\item
$G'\in \mathcal{G}_{m,(q-1)}$ with symptoms $B_{m,(q-1)}=B\backslash\{b_q\}$ such that $pa_{G'}(b_i)=pa_{G}(b_i)$, $\forall\ b_i\in B_{m,(q-1)}$. This implies $c_{G'}\in \mathbf{P}_{m,q-1}$;
\item
$G'' \in \mathcal{G}_{m,1}$ with symptom $B_{m,1}=\{b_q\}$ such that $pa_{G''}(b_q)=pa_{G}(b_q)$. This implies $c_{G''}\in \mathbf{P}_{m,1}$.
\end{itemize}
With a proper permutation of coordinates, we can write $c_G$ in the form of:
\begin{center}
$c_G=(c_{G'}, c_{G''})$.
\end{center}
Recall that $vert(\mathbf{P}_{m,q}) = \{c_G: G\in \mathcal{G}_{m,q}\}$, so $\forall x \in \mathbf{P}_{m,q}$, with the same permutation of coordinates, we have:
\begin{equation}\label{eq:diagprod1}
x=\sum\limits_{G\in\mathcal{G}_{m,q}} \alpha_G c_G =
 (\sum\limits_{G\in\mathcal{G}_{m,q}} \alpha_G c_{G'} \ ,\  \sum\limits_{G\in\mathcal{G}_{m,q}} \alpha_G c_{G''})
\end{equation}
where $0 \leq \alpha_G\leq 1$, $\forall G\in \mathcal{G}_{m,q}$ and $\sum\limits_{G\in\mathcal{G}_{m,q}} \alpha_G =1$.

Note that $\sum\limits_{G\in\mathcal{G}_{m,q}} \alpha_G c_{G'} \in \mathbf{P}_{m,q-1}$ and $\sum\limits_{G\in\mathcal{G}_{m,q}} \alpha_G c_{G''}\in \mathbf{P}_{m,1}$, Equation \eqref{eq:diagprod1} implies $x\in \mathbf{P}_{m,q-1}\times \mathbf{P}_{m,1}$.
Hence:
\[
\mathbf{P}_{m,q}\subseteq \mathbf{P}_{m,q-1}\times \mathbf{P}_{m,1}.
\]

Second, we need to prove: $\mathbf{P}_{m,q-1}\times \mathbf{P}_{m,1} \subseteq \mathbf{P}_{m,q}$.

Let $\mathcal{G}_{m,q-1}$ has symptoms $B_{m,(q-1)}=B\backslash\{b_q\}$ and $\mathcal{G}_{m,1}$ has symptom $B_{m,1}=\{b_q\}$. $\forall \ G'\in \mathcal{G}_{m,(q-1)}$ and $G'' \in \mathcal{G}_{m,1}$, we can define $G\in \mathcal{G}_{m,q}$ such that $pa_{G}(b_i)=pa_{G'}(b_i)$, $\forall\ b_i\in B_{m,(q-1)}$, and $pa_{G}(b_q)=pa_{G''}(b_q)$. $c_G$ has the form of $c_G=(c_{G'}, c_{G''})$.

$\forall\ x\in \mathbf{P}_{m,q-1}\times \mathbf{P}_{m,1}$, $x$ can be written as:
\begin{center}
$\begin{array}{ccl}
x & = &
 (\sum\limits_{G'\in\mathcal{G}_{m,q-1}} \beta_{G'} c_{G'} \ ,\  \sum\limits_{G''\in\mathcal{G}_{m,1}} \gamma_{G''} c_{G''})
=\sum\limits_{G'\in\mathcal{G}_{m,q-1}} \sum\limits_{G''\in\mathcal{G}_{m,1}}
   \beta_{G'} \gamma_{G''} (c_{G'}\ ,\  c_{G''}) \\
   & = & \sum\limits_{G'\in\mathcal{G}_{m,q-1}} \sum\limits_{G''\in\mathcal{G}_{m,1}}
  ( \beta_{G'} \gamma_{G''}) \  c_G\ ,
\end{array}$
\end{center}
where $0 \leq \beta_{G'}, \gamma_{G''}\leq 1$, $\forall G'\in \mathcal{G}_{m,q-1}$, $\forall G''\in \mathcal{G}_{m,1}$, and $\sum_{G'\in\mathcal{G}_{m,q-1}} \beta_{G'} =1$, $\sum_{G''\in\mathcal{G}_{m,1}} \gamma_{G''} =1$.
Note that
\[
\sum\limits_{G'\in\mathcal{G}_{m,q-1}} \sum\limits_{G''\in\mathcal{G}_{m,1}}
  ( \beta_{G'} \gamma_{G''}) = \sum\limits_{G'\in\mathcal{G}_{m,q-1}} \beta_{G'}( \sum\limits_{G''\in\mathcal{G}_{m,1}}
   \gamma_{G''})= \sum\limits_{G'\in\mathcal{G}_{m,q-1}}\beta_{G'}=1,
\]
which leads to $x\in \mathbf{P}_{m,q}$. Hence:
\[
\mathbf{P}_{m,q-1}\times \mathbf{P}_{m,1} \subseteq \mathbf{P}_{m,q}.
\]
Therefore,
\[
\mathbf{P}_{m,q} = \mathbf{P}_{m,q-1}\times \mathbf{P}_{m,1}= \underbrace{\Delta_{2^m-1}  \times \cdots \times \Delta_{2^m-1}}_{\mbox{q-1 many}}\times \Delta_{2^m-1} = \underbrace{\Delta_{2^m-1}  \times \cdots \times \Delta_{2^m-1}}_{\mbox{q many}}.
\]

\end{itemize}
\end{proof}

Theorem \ref{thm:diagprod} implies that $\mathbf P_{m,n}$ is a simple polytope with dimension $n\cdot (2^m-1)$. In Section \ref{proof:simple}, we will give another proof which use linear algebra to show that $\mathbf P_{m,n}$ is simple and obtain its dimension.
(cim-polytope)

\subsection{Expression of facets of $\mathbf{P}_{m,n}$}

Based on Theorem \ref{thm:diagprod}, we are going to show the expression of facets of $\mathbf{P}_{m,n}$ using the following lemma:

\begin{lem} \label{lem:diagprod} \cite{Ziegler}
Suppose $\mathbf{P}$ is the direct product of simplices $\Delta_{\alpha_1},\ldots,\Delta_{\alpha_k}$. Then every facet of $\mathbf{P}$ has the form of $\Delta_{\alpha_1}\times \ldots \times \Delta_{\alpha_{i-1}} \times F_{\alpha_i}\times \Delta_{\alpha_{i+1}}\times \ldots \times \Delta_{\alpha_k}$, where $F_{\alpha_i}$ is a facet of $\Delta_{\alpha_i}$.
\end{lem}

\begin{rem}
Lemma \ref{lem:diagprod} implies that in order to study the facets of a direct product of simplices, we can simply study the facets of each simplex. As by Theorem \ref{thm:diagprod}, $\mathbf{P}_{m,n}$ is a direct product of $n$ many $\mathbf{P}_{m,1}$, our problem is simplified as studying the facets of $\mathbf{P}_{m,1}$. Thus we assume $B=\{b_1\}$ in the following content of this section.
\end{rem}

Assume $A=\{a_1, \ldots, a_m\}$ and $B=\{b_1\}$. 
By Proposition \ref{prop:diag4}, the vertices of $\mathbf{P}_{m,1}$ has at most $2^m-1$ many non-zero coordinates. We define the indeterminates, i.e. variables, $\{x_s,\ s\subseteq A,\ s\neq \emptyset\}$, where one indeterminate $x_s$ for each coordinate $c_G(s \cup \{b_1\})$ in the characteristic imset $c_G$, $G \in \mathcal G_{m,1}$. Define the vector of indeterminates $x=\{x_s,\ s\subseteq A,\ s\neq \emptyset\}$. Suppose $A_m x \leq b_m$ is the system of inequalities that defines $\mathbf{P}_{m,1}$. We can define a $2^m\times 2^m$ matrix: $D_m$ = $[b_m|-A_m]$. Denote the elements in $D_m$ by $(d_{st})_{s\subseteq A, t\subseteq A}$ so that we can rewrite the system of inequalities as: $d_{s\emptyset}+\sum_{t\subseteq A, t\neq \emptyset} d_{st}x_t \geq 0,\ s\subseteq A$. Then we have the expression of $2^m$ facets of $\mathbf{P}_{m,1}$ as following:
\[
F_s = \mathbf{P}_{m,1}\cap \{x: d_{s\emptyset}+\sum_{t\subseteq A, t\neq \emptyset} d_{st}x_t = 0\},\ s\subseteq A,
\]
where the elements $d_{st}$, $s$, $t\subseteq A$ can be obtained using Theorem \ref{thm:diagfacet}.

\begin{thm}\label{thm:diagfacet}
The elements in matrix $D_m$ satisfies:
\begin{itemize}
\item
$d_{st} \neq 0$ if and only if $s\subseteq t$;
\item
if $s\subseteq t$, then $d_{st}=(-1)^{|t|-|s|}$.
\end{itemize}
This implies that $\mathbf{P}_{m,1}$ has $2^m$ facets:
\[
F_s = \mathbf{P}_{m,1}\cap \{x: d_{s\emptyset}+\sum_{t\subseteq A, t\neq \emptyset} d_{st}x_t = 0\},\ s\subseteq A.
\]
What's more, $\forall \ s\subseteq A$, $vert(\mathbf{P}_{m,1})\backslash \{c_{G_s}\} \subset F_s$, where $pa_{G_s}(b_1)=s$.
\end{thm}

\begin{proof}
For convenience, let $x_\emptyset \equiv 1$. $\forall s\subseteq A$, let $d_{s\cdot} = (d_{st})_{t \subseteq A}$ be the corresponding row of $D_m$, and $G_s$ be the graph in $\mathcal G_{m,1}$ such that $pa_{G_s}(b_1)=s$. Now we can rewrite the system of inequalities as: 
\[
\sum\limits_{t\subseteq A} d_{st}x_t=d_{s\cdot}(1\ x)^T \geq 0, \mbox{ for }\forall\ s\subseteq A.
\]
We are going to prove that $\forall s\subseteq A$, we can find $2^m-1$ vertices on $F_s$ that are linearly independent, and this implies that $F_s$ is a facet of $\mathbf{P}_{m,1}$. In fact, we will prove that: $\{c_{G_{s'}}, \ s'\subseteq A, \ s'\neq s\}\subset F_s$ and $c_{G_s}\notin F_s$, i.e. $d_{s\cdot}(1\ c_{G_{s'}})^T = 0,\ \forall \ s'\subseteq A,\ s'\neq s$ and $d_{s\cdot}(1\ c_{G_s})^T > 0$.\\
Notice that $\forall t\subseteq A$, $c_{G_{s'}}(t\cup \{b_1\})\neq 0$ if and only if $t\subseteq pa_{c_{G_{s'}}}(b_1)=s'$, and $d_{st} \neq 0$ if and only if $s\subseteq t$. So:
\[
d_{s\cdot}(1\ c_{G_{s'}})^T=d_{s\emptyset}+\sum\limits_{t\subseteq A,\ t\neq \emptyset} d_{st} c_{G_{s'}}(t\cup\{b_1\}) = d_{s\emptyset}+\sum\limits_{s\subseteq t\subseteq s',\ t\neq \emptyset} d_{st} = \sum\limits_{s\subseteq t\subseteq s'} d_{st} .
\]
Therefore, we have: 
\begin{itemize}
\item
if $s=s'$, then $d_{s\cdot}(1\ c_{G_{s'}})^T=d_{ss}=1>0$;
\item
if $s\subsetneq s'$, then $d_{s\cdot}(1\ c_{G_{s'}})^T=\sum\limits_{s\subseteq t\subseteq s'} (-1)^{|t|-|s|}=\sum\limits_{t'\subseteq s'\backslash s} (-1)^{|t'|}=0$;
\item
if $s \nsubseteq s'$, then $d_{s\cdot}(1\ c_{G_{s'}})^T=0$.
\end{itemize}
\end{proof}

\begin{ex}[Facets of $\mathbf P_{2,1}$]
Notation adopted from Theorem \ref{thm:diagfacet}. Fix $m=2$ and $n=1$.

\begin{multicols}{2}
\begin{center}
\scalebox{0.8}{
\includegraphics{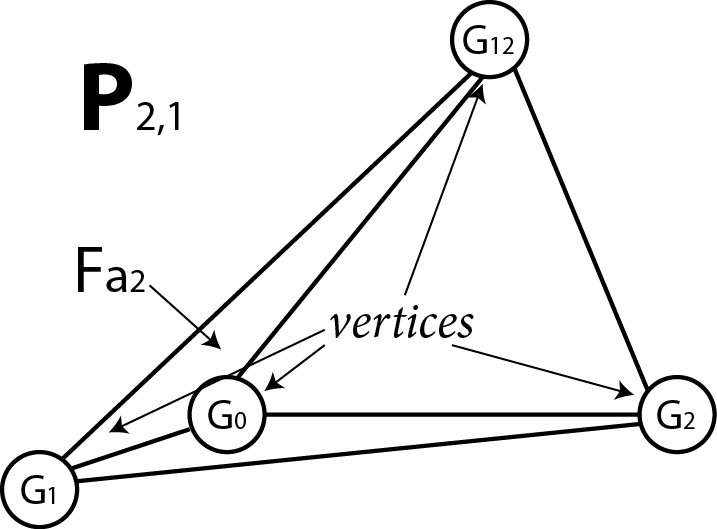}
}
\end{center}
All characteristic imsets are given as a matrix:
$$ \bordermatrix{ & \cr
	&  c_{G_0} \cr	&  c_{G_1} \cr	&  c_{G_2} \cr	&  c_{G_{12}} } 
= \bordermatrix{ T & a_1b_1 & a_2b_1 & a_1a_2b_1\cr
	       &  0  &  0  &  0   \cr
                &  1  &  0  &  0   \cr
	       &  0  &  1  &  0   \cr
	       &  1  &  1  &  1  
}$$
The matrix $D_2=[b_2 | -A_2]$:
$$D_2=  \bordermatrix{  s\backslash t  & \emptyset  & a_1 & a_2 & a_1a_2\cr
	     \emptyset   &  1  &  -1  &  -1 &   1   \cr
              a_1              &  0  &   1   &   0  & -1   \cr
	     a_2              &  0  &   0   &   1  & -1   \cr
	     a_1a_2        &  0  &   0   &   0  &   1
}$$
The system of inequalities that defines $\mathbf{P}_{2,1}$:
$$\bordermatrix{  s\backslash t  & \emptyset  & a_1 & a_2 & a_1a_2 & \cr
	     \emptyset   &  1  &  -x_{a_1}  &  -x_{a_2} &   +x_{a_1a_2}  & \geq 0  \cr
              a_1              &     &   x_{a_1}   &      & -x_{a_1a_2}  & \geq 0 \cr
	     a_2              &     &       &   x_{a_2}  & -x_{a_1a_2}  & \geq 0 \cr
	     a_1a_2        &     &       &      &   x_{a_1a_2}  & \geq 0
}$$
\end{multicols}
Vertices $c_{G_0}$, $c_{G_1}$ and $c_{G_{12}}$ are in the facet $F_{a_2}$ while $c_{G_2}$ is not.
\end{ex}

\section{The cim-polytopes for Bayesian networks}\label{BN:generalize}

The results in Section \ref{BN:diagnosis} are limited to diagnosis models. In this section, we will generalize the results to all Bayesian networks with the same underlying order.

\subsection{Underlying ordering of DAGs}

For a set of random variables $N=\{a_1,\ldots,a_n\}$, where now $n$ is
the total number of nodes in $N$. $\forall G\in DAGs(N)$, there exists
an {\em underlying ordering} over $N$, $[n]_G=(a_{[1]}, \ldots, a_{[n]})$, such that if $a_{[i]} \rightarrow a_{[j]}$ in $G$, then $i<j$.  We are are now interested in the class of graphs which share a specific underlying ordering $[n]$, i.e. $\mathcal{G}_{[n]}=\{G\in DAGs(N): [n]_G=[n]\}$, and its cim-polytope $\mathbf{P}_{[n]}=\mathbf P_{\mathcal{G}_{[n]}}$.

\begin{ex} [Underlying ordering of graphs]
Let $N=\{a_1,a_2,a_3\}$. Consider an ordering over $N$, $[n]=(a_2, a_1, a_3)$, i.e. $a_{[1]} = a_2$, $a_{[2]}=a_1$ and $a_{[3]}=a_3$. Then $\forall G \in \mathcal{G}_{[n]}$, the only type of directed edges allowed in $G$ are $a_{[i]} \rightarrow a_{[j]}$, where $i<j$. For instance, $a_2 \rightarrow a_1$ is allowed while $a_1 \rightarrow a_2$ is not. Thus graph $G_1$ in Figure \ref{fig:order:G1} and graph $G_2$ in Figure \ref{fig:order:G2} are both in $\mathcal{G}_{[n]}$. Graph $G_3$ in Figure \ref{fig:order:G3} is not in $\mathcal{G}_{[n]}$ since it has arrow $a_1 \rightarrow a_2$, and the underlying ordering for $G_3$, i.e. $[n]_{G_3}$, can either be $(a_1, a_2, a_3)$ or $(a_1, a_3, a_2)$.
\begin{figure}
  \centering
  \subfigure[$G_1$]{
    \label{fig:order:G1} 
    \includegraphics[width=1.25in]{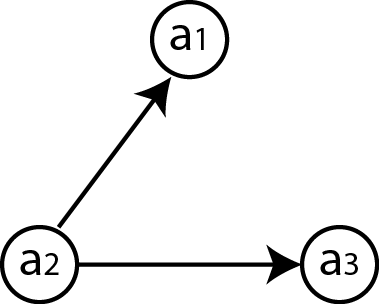}}
  \hspace{0.5in}
  \subfigure[$G_2$]{
    \label{fig:order:G2} 
    \includegraphics[width=1.25in]{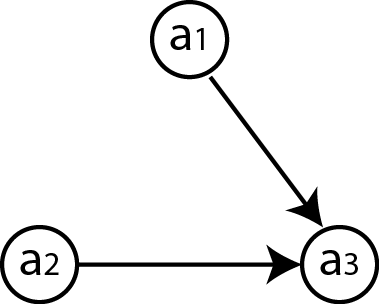}}
  \hspace{0.5in}
  \subfigure[$G_3$]{
    \label{fig:order:G3} 
    \includegraphics[width=1.25in]{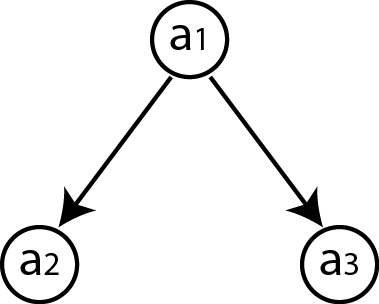}}
  \caption{Three graphs to illustrate the underlying ordering of graphs}
  \label{fig:order} 
\end{figure}

\end{ex}

\begin{rem}\label{rem:gen1}
For a specific ordering $[n]$ and an arbitrary $G \in \mathcal{G}_{[n]}$, we have the following proposition that is similar with Proposition \ref{prop:diag2}.
\begin{itemize}
\item
$\forall T\subseteq N$, $|T|=k \geq 2$, we can order the elements in $T$ according to $[n]$ and write $T$ in the form of $a_{[i_1]} a_{[i_2]} \ldots a_{[i_k]}$ where $i_1 < i_2 < \cdots <i_k$. Then $c_G(T)=\prod_{s=i_1,\ldots,i_{k-1}} c_G(a_{[s]} a_{[i_k]})$. This property means that the whole $c_G$ is determined by $n \choose 2$ coordinates, $\{c_G(a_{[i]} a_{[j]})$, $i<j\}$, which can also be interpreted as the existence of the directed edges $a_{[i]} \rightarrow a_{[j]}$, $i < j$.
\end{itemize}
Another way to see this property is that $\forall\ G \in \mathcal{G}_{[n]}$, $G$ can be determined by $pa_G(a_{[i]})$, $i = 2, \ldots, n$ since $pa_G(a_{[1]})=\emptyset$. Similarly with Remark \ref{rem:diag1}, we can consider a permutation of coordinates in $c_G$ that corresponds to a permutation of $T$, then these coordinates can be broken into $n-1$ parts:
\[
\underline{(12)}, \uwave{(13), (23), (123)}, \underline{(14), (24), (34),\ldots, (1234)}, \ldots, \underline{(1n), (2n), \ldots, ((n-1)n), \ldots, (12\ldots n)}
\]
where $(i_1 \ldots i_k)$ stands for $T = a_{[i_1]} a_{[i_2]} \ldots a_{[i_k]} $, $\{i_1,\ldots,i_k\} \subseteq \{1, \ldots, n\}$. The $k$-th part of the coordinations, $\{c_G(T)$: $a_{[j]} \notin T$, $\forall j>k\}$ only depend on $pa_G(a_{[k]})$, and different parts are completely irrelevant in the sense that $pa_G(a_{[k]})$, $a_{[k]} \in N$, can be decided separately.

\end{rem}

\subsection{Structure, edges and facets of $\mathbf P_{[n]}$}

\begin{thm}\label{thm:genprod}
Suppose $n \geq 2$. $\mathbf{P}_{[n]}$ is a direct product of a sequence of simplices:
\begin{center}
$\mathbf{P}_{[n]}= \underbrace{ \Delta_{2^1-1} \times \Delta_{2^2-1} \times \cdots \times \Delta_{2^{n-1}-1}}_{\mbox{n-1 simplices}}$,
\end{center}
where the $i_{th}$ simplex $\Delta_{2^{i}-1}$ is the same with the cim-polytope for diagnosis models, $\mathbf{P}_{i,1}$, with diseases $A=\{a_{[1]}, \ldots, a_{[i]}\}$ and one symptom $\{a_{[i+1]}\}$.
\end{thm}

\begin{proof}
See Section \ref{proof:generalize}.
\end{proof}

\begin{rem}
Two immediate results from Theorem \ref{thm:genprod} are:
\begin{itemize}
\item
the dimension of $\mathbf{P}_{[n]}$ is $2^n - (n+1)$, and it is a simple polytope;
\item
the facets of $\mathbf{P}_{[n]}$ can be obtained by Lemma \ref{lem:diagprod} and Theorem \ref{thm:diagfacet}.
\end{itemize}
\end{rem}

\begin{rem}
Note that the equality in Theorem \ref{thm:genprod} is actually $\mathbf{P}_{[n]}=\Delta_{2^0-1}\times \Delta_{2^1-1} \times \Delta_{2^2-1} \times \cdots \times \Delta_{2^{n-1}-1}$, where $\Delta_{2^0-1}$ is omitted as it has dimension 0 (a point).
Theorem \ref{thm:genprod} and its proof also imply that $\forall x \in \mathbf{P}_{[n]}$, $x \in vert(\mathbf{P}_{[n]})$ if and only if with the permutation of coordinates in Remark \ref{rem:gen1}, $x$ can be written in the form of $x=(v_1, v_2,\ldots,v_{n-1})$, where $v_i$ is the vertex of $\Delta_{2^{i}-1}$, $i=1,\ldots,n-1$. Suppose $x=c_G$, $G\in \mathcal{G}_{[n]}$, then $v_i=c_{G_i}$, where $G_i$ is in $\mathcal{G}_{i,1}$ with diseases $N_{[i]}$ and symptom $a_{[i+1]}$, $i=1,\ldots,n-1$, and $pa_{G_i}(a_{[i+1]})=pa_{G}(a_{[i+1]})$.
\end{rem}

 The following theorem will be stated in two forms which are equivalent by Theorem \ref{thm:genprod} and Lemma \ref{lem:onesymp}.
\begin{thm}\label{thm:genedge}
Fix an underlying ordering $[n]$ over $N$. 
\begin{itemize}
\item
(From the view of graph theory.)
Two graphs, $G_1$, $G_2\in \mathcal{G}_{[n]}$ are neighbors in $\mathcal{G}_{[n]}$ if and only if: $\exists\ a_{[i]}\in N$ such that $pa_{G_1}(a_{[i]})\neq pa_{G_2}(a_{[i]})$ and $pa_{G_1}(a_{[j]}) = pa_{G_2}(a_{[j]})$, $\forall\ a_{[j]}\in N$ and $a_{[j]}\neq a_{[i]}$, i.e., all nodes but one have exactly the same parent sets in both $G_1$ and $G_2$.
\item
(From the view of polyhedral geometry.)
$\forall \ \mathbf x \in \mathbf{P}_{[n]}$, $\mathbf x$ is on an edge of $\mathbf{P}_{[n]}$ if and only if with the permutation of coordinates showed in Remark \ref{rem:gen1} $\mathbf x$ can be written in the form of $\mathbf x=(v_1,\ldots,v_{i-1},e_i,v_{i+1},\ldots,v_{n-1})$, where $e_i$ belongs to an edge on $\Delta_{2^i-1}$, $i\in\{1,\ldots, n-1\}$, and $v_j \in vert(\Delta_{2^j-1})$, $j\in\{1,\ldots, n-1\}\backslash \{i\}$.
\end{itemize}
\end{thm}

\begin{proof}
See Section \ref{proof:generalize}.
\end{proof}

\subsection{Graphes with forbidden (or fixed) edges}

Fix an underlying ordering of nodes $[n]$ and consider $\mathcal G_{[n]}$. When a specific set of directed edges are forbidden in $\mathcal G_{[n]}$, we can define sets of nodes $\Omega = \{\Omega_i ^0$, $i=2, \ldots, n\} \cup \{\Omega_i ^1$, $i=2, \ldots, n\}$ such that $\Omega_i^0 \subseteq \Omega_i^1 \subseteq \{a_{[1]}, \ldots, a_{[i-1]}\}$, and the class of graphs we are interested in becomes $\mathcal G_{[n], \Omega}=\{G\in DAGs(N): [n]_G = [n]$, $\Omega_i^0 \subseteq pa_G(a_{[i]}) \subseteq \Omega_i^1$, $i=2, \ldots, n \}$, i.e. edges $\{ a_{[j]}\rightarrow a_{[i]}$: $a_{[j]}\in \Omega_i^0$, $i=2, \ldots, n\}$ are fixed edges, and edges $\{ a_{[j]}\rightarrow a_{[i]}$: $a_{[j]}\in \{ a_{[1]}, \ldots, a_{[i-1]} \}\backslash \Omega_i^1$, $i=2, \ldots, n\}$ are forbidden edges. The cim-polytope for $\mathcal G_{[n], \Omega}$ is $\mathbf P_{\mathcal G_{[n], \Omega}}$. Using similar strategy, we are able to show that $\mathbf P_{\mathcal G_{[n], \Omega}}$ is a direct product of a sequence of simplices:

\begin{thm}
\begin{align}\label{eq:genprod2}
& \mathbf P_{\mathcal G_{[n], \Omega}} & = &  \mathbf P_{a_{[2]}} \times \ldots \times\mathbf P_{a_{[n]}} \nonumber \\
& & = & \underbrace{\Delta_{2^{|\Omega_2^1|-|\Omega_2^0|}-1} \times \cdots \times \Delta_{2^{|\Omega_2^1|-|\Omega_2^0|}-1}}_{2^{|\Omega_2^0|} \mbox{ many}} \times \ldots \times  \underbrace{\Delta_{2^{|\Omega_n^1|-|\Omega_n^0|}-1} \times \cdots \times \Delta_{2^{|\Omega_n^1|-|\Omega_n^0|}-1}}_{2^{|\Omega_n^0|} \mbox{ many}} ,
\end{align}
where the $i$-th polytope $\mathbf P_{a_{[i+1]}} = \underbrace{\Delta_{2^{|\Omega_{i+1}^1|-|\Omega_{i+1}^0|}-1} \times \cdots \times \Delta_{2^{|\Omega_{i+1}^1|-|\Omega_{i+1}^0|}-1}}_{2^{|\Omega_{i+1}^0|} \mbox{ many}} $ is a $(2^{|\Omega_{i+1}^1|}-2^{|\Omega_{i+1}^0|})$-face of $\mathbf P_{|\Omega_{i+1}^1|,1}=\Delta_{2^{|\Omega_{i+1}^1|}-1}$, where $\mathbf P_{|\Omega_{i+1}^1|,1}$ is the cim-polytope for diagnosis models with diseases $A = \Omega_{i+1}^1$ and one symptom $a_{[i+1]}$.

\end{thm}

\begin{proof}
To prove $\mathbf P_{a_{[i+1]}} = \underbrace{\Delta_{2^{|\Omega_{i+1}^1|-|\Omega_{i+1}^0|}-1} \times \cdots \times \Delta_{2^{|\Omega_{i+1}^1|-|\Omega_{i+1}^0|}-1}}_{2^{|\Omega_{i+1}^0|} \mbox{ many}} $, we permutate the coordinates in the following way:

\begin{align}
\{T: T \subseteq \Omega_{i+1}^0 \cup a_{[i+1]}\} \cup \bigcup_{\Omega_{s} \subseteq \Omega_{i+1}^0} \{T\subseteq \Omega_{i+1}^1 \cup a_{[i+1]}: T\cap  \Omega_{i+1}^0 = \Omega_{s} \},
\end{align}
i.e. $c_G(T)$, $\forall G \in \mathcal G_{[n],\Omega}$, can be split into the following subvectors: $(c_G(T)$, where $T \subseteq \Omega_{i+1}^0 \cup a_{[i+1]})$, $(c_G(T)$, where $T\subseteq \Omega_{i+1}^1 \cup a_{[i+1]}$ and $T\cap  \Omega_{i+1}^0 = \Omega_{s})$, $\forall \Omega_{s} \subseteq \Omega_{i+1}^0$.

Then use the strategy similar with the previous proofs, we can prove the following:
\begin{itemize}
\item
$c_G(T)$, $T\subseteq \Omega_{i+1}^0 \cup a_{[i+1]}$, are all fixed;
\item
$\forall \Omega_{s} \subseteq \Omega_{i+1}^0$, the convex hull of $\{(c_G(T)$, where $T\subseteq \Omega_{i+1}^1 \cup a_{[i+1]}$ and $T\cap  \Omega_{i+1}^0 = \Omega_{s})$: $\forall G \in \mathcal G_{[n],\Omega} \}$ is $\Delta_{2^{|\Omega_{i+1}^1|-|\Omega_{i+1}^0|}-1}$ (see Example \ref{ex:fixedge});
\item
$\mathbf P_{a_{[i+1]}} = \underbrace{\Delta_{2^{|\Omega_{i+1}^1|-|\Omega_{i+1}^0|}-1} \times \cdots \times \Delta_{2^{|\Omega_{i+1}^1|-|\Omega_{i+1}^0|}-1}}_{2^{|\Omega_{i+1}^0|} \mbox{ many}} $;
\item
Equation \ref{eq:genprod2} holds.
\end{itemize}

\end{proof}

\begin{ex}\label{ex:fixedge}
Consider a DAG $G$ which has $7$ nodes $\{a_1, \ldots, a_7\}$. After fix an underlying ordering, we can write these nodes as $\{a_{[1]}, \ldots, a_{[7]}\}$, where $a_{[i]} \rightarrow a_{[j]}$ in $G$ implies $i<j$. Suppose edges $a_{[1]} \rightarrow a_{[6]}$ and $a_{[2]} \rightarrow a_{[6]}$ are fixed and edge $a_{[5]} \rightarrow a_{[6]}$ is forbidden. Then coordinates $c_G(T)=0$ if $a_{[5]} \in T$, and other coordinates $c_G(T)$ where $a_{[j]} \notin T$, $\forall j > 6$, can be ordered as following (values with respect to different DAGs are listed as a matrix):
\tiny{
\[
 \begin{array}{c}  T\backslash \{a_{[6]}\}   \\ \\ \\ 
\end{array}
\begin{bmatrix}{ \begin{smallmatrix}  a_{[1]} & a_{[2]} & a_{[1]} a_{[2]} & a_{[3]} & a_{[4]} & a_{[3]} a_{[4]} & a_{[1]} a_{[3]} & a_{[1]} a_{[4]} & a_{[1]} a_{[3]} a_{[4]} & a_{[2]} a_{[3]} & a_{[2]} a_{[4]} & a_{[2]} a_{[3]} a_{[4]} & a_{[1]} a_{[2]} a_{[3]} & a_{[1]} a_{[2]} a_{[4]} & a_{[1]} a_{[2]} a_{[3]} a_{[4]} \cr
	         1  &  1  &  1 & 0 & 0 & 0 & 0 & 0 & 0 & 0 & 0 & 0 & 0 & 0 & 0 \cr
	         1  &  1  &  1 & 1 & 0 & 0 & 1 & 0 & 0 & 1 & 0 & 0 & 1 & 0 & 0 \cr
	         1  &  1  &  1 & 0 & 1 & 0 & 0 & 1 & 0 & 0 & 1 & 0 & 0 & 1 & 0 \cr
	         1  &  1  &  1 & 1 & 1 & 1 & 1 & 1 & 1 & 1 & 1 & 1 & 1 & 1 & 1
\end{smallmatrix}}
\end{bmatrix}
\]
}\normalsize
where the $4$ rows correspond to graphs $G_i$, $i = 1, \ldots, 4$, such that $pa_{G_1}(a_{[6]})=\{a_{[1]}, a_{[2]}\}$, $pa_{G_1}(a_{[6]})=\{a_{[1]}, a_{[2]}, a_{[3]}\}$, $pa_{G_1}(a_{[6]})=\{a_{[1]}, a_{[2]}, a_{[4]}\}$ and $pa_{G_1}(a_{[6]})=\{a_{[1]}, a_{[2]}, a_{[3]}, a_{[4]}\}$.

\end{ex}

It is obvious that the cim-polytope for diagnosis models, $\mathbf P_{m,n}$, is a special case of $\mathbf P_{\mathcal G_{[n], \Omega}}$: the underlying ordering of nodes is $(a_1, \ldots, a_m, b_1, \ldots, b_n)$ (the ordering is not unique in the sense that the order of two diseases or two symptoms can exchange), $\Omega_i^0 = \Omega_i^1  = \emptyset$ for $i=1, \ldots, m$, while $\Omega_i^0 = \emptyset$ and $\Omega_i^1 = \{a_1, \ldots, a_m\}$ for $i = m+1, \ldots, m+n$. Note that based on Equation \eqref{eq:genprod2}, all edges of $\mathbf P_{\mathcal G_{[n], \Omega}} $ can be found similarly with Theorem \ref{thm:genedge}, and the its facets can also be obtained by Lemma \ref{lem:diagprod} and Theorem \ref{thm:diagfacet}.



\section{Proofs in Section \ref{prop}}\label{proof:diagnosis}

\subsection{Proof of Lemma \ref{lem:onesymp}}

\begin{proof}
Let $N=A\cup B$, where $A=\{a_1, \ldots, a_m\}$ and $B=\{ b_1\}$. We need to prove: $\exists$ a cost vector $w$, such that $w\cdot c_{G_1} = w\cdot c_{G_2} > w\cdot c_{G_3}$, $\forall\ G_3\in \mathcal{G}_{m,1}$ distinct with $G_1$ and $G_2$.\\
By Remark \ref{rem:diag1}, $G_1$ and $G_2$ are determined by $pa_{G_1}(b_1)$ and $pa_{G_2}(b_1)$, respectively. We will discuss by two scenarios of $pa_{G_1}(b_1)$ and $pa_{G_2}(b_1)$: one is a subset of the other, and neither one is a subset of the other.
\begin{itemize}
\item[(1)]
One is a subset of the other. WLOG, suppose $pa_{G_1}(b_1)\subsetneq pa_{G_2}(b_1)$.

Define: $A_1=pa_{G_1}(b_1)$, $A_2=pa_{G_2}(b_1)$, $A_{2\backslash 1}=pa_{G_2}(b_1)\backslash pa_{G_1}(b_1)$, and $A_{comp}=(pa_{G_2}(b_1))^c$ (i.e. the complement set of $pa_{G_2}(b_1)$). Note that: $A_{2\backslash 1} \neq \emptyset$, $A_1$ and $A_{comp}$ can be $\emptyset$; $A_1$, $A_{2\backslash 1}$ and $A_{comp}$ is a partition of $N$.

Consider a function $w: \mathcal P(N) \mapsto \mathbb R$ where $w(T)=0$ if $|T| <2$. Then similar with imsets, $w$ can also be considered as a vector, and we assume that the permutations of coordinates in $w$ and in characteristic imsets coincide. 
\begin{itemize}
\item
If $|A_{2\backslash 1}|>1$, we define $w$ as:
\[
w(T) = \left\{ \begin{array}{lcl}
c   & \mbox{for} & T=a_i b_j,\ a_i\in A_1 \\
-c & \mbox{for} & T=a_i b_j,\ a_i\notin A_1 \\
|A_{2\backslash 1}|\cdot c & \mbox{for} & T=A_{2\backslash 1}\cup\{b_1\} \\
0   & \mbox{for} & T\subset N, |T|>2,\mbox{ and } T\neq A_{2\backslash 1}\cup\{b_1\} 
\end{array}\right.
\]
where $c$ is a positive number.\\
Then $\forall\ G_3 \in \mathcal{G}_{m,1}$, we have:
\small{\[
\begin{array}{ccl}
w\cdot c_{G_3} & = & |A_1\cap pa_{G_3}(b_1)|\cdot c-|pa_{G_3}(b_1)\backslash A_1|\cdot c
+|A_{2\backslash 1}|\cdot c \cdot c_{G_3}(A_{2\backslash 1}\cup\{b_1\})\\
 & = & |A_1\cap pa_{G_3}(b_1)|\cdot c-|pa_{G_3}(b_1)\cap A_{2\backslash 1}|\cdot c -|pa_{G_3}(b_1)\cap A_{comp}|\cdot c \\
&   & +|A_{2\backslash 1}|\cdot c \cdot c_{G_3}(A_{2\backslash 1}\cup\{b_1\}).
\end{array}
\]}
\normalsize
In this equation:
\begin{itemize}
\item
$ |A_1\cap pa_{G_3}(b_1)|\cdot c \leq |A_1|\cdot c$, where ``=" holds if and only if $A_1\subset pa_{G_3}(b_1)$;
\item
$-|pa_{G_3}(b_1)\cap A_{2\backslash 1}|\cdot c+|A_{2\backslash 1}|\cdot c \cdot c_{G_3}(A_{2\backslash 1}\cup\{b_1\})\leq 0$, where ``=" holds if and only if $pa_{G_3}(b_1)\cap A_{2\backslash 1} = \emptyset$ or $A_{2\backslash 1}$;
\item
$-|pa_{G_3}(b_1)\cap A_{comp}|\cdot c \leq 0$, where ``=" holds if and only if $pa_{G_3}(b_1)\cap A_{comp}=\emptyset$.
\end{itemize}
Therefore, $w\cdot c_{G_3}\leq |A_1|\cdot c$, where ``=" holds if and only if $G_3=G_1$ or $G_2$.

\item
If $|A_{2\backslash 1}|=1$, we let $A_{2\backslash 1}=\{a_q\}$, and define $w$ as:
\[
w(T) = \left\{ \begin{array}{lcl}
c   & \mbox{for} & T=a_i b_j,\ a_i\in A_1 \\
-c & \mbox{for} & T=a_i b_j,\ a_i\notin A_2 \\
0   & \mbox{for} & T=a_q b_1 \\
0   & \mbox{for} & T\subset N, |T|>2, \mbox{ and } T\neq A_{2\backslash 1}\cup\{b_1\} 
\end{array}\right.
\]
where $c$ is a positive number.

Then $\forall\ G_3 \in \mathcal{G}_{m,1}$, we have:
\[
w\cdot c_{G_3} = |A_1\cap pa_{G_3}(b_1)|\cdot c-|pa_{G_3}(b_1)\cap A_{comp}|\cdot c.
\]
Again, in this equation:
\begin{itemize}
\item
$ |A_1\cap pa_{G_3}(b_1)|\cdot c \leq |A_1|\cdot c$, where ``=" holds if and only if $A_1\subset pa_{G_3}(b_1)$;
\item
$-|pa_{G_3}(b_1)\cap A_{comp}|\cdot c \leq 0$, where ``=" holds if and only if $pa_{G_3}(b_1)\cap A_{comp}=\emptyset$.
\end{itemize}
To satisfy the above two conditions, we must have $pa_{G_3}(b_1)=A_1$ or $(A_1\cup {a_q})$. Therefore, again, we have: $w\cdot c_{G_3}\leq |A_1|\cdot c$, where ``=" holds if and only if $G_3=G_1$ or $G_2$.
\end{itemize}

\item[(2)]
Neither one is a subset of the other.

Define: $A_1=pa_{G_1}(b_1)$, $A_2=pa_{G_2}(b_1)$, $A_{1\cap 2}=pa_{G_1}(b_1)\cap pa_{G_2}(b_1)$, $A_{1\backslash 2}=pa_{G_1}(b_1)\backslash pa_{G_2}(b_1)$, $A_{2\backslash 1}=pa_{G_2}(b_1)\backslash pa_{G_1}(b_1)$, $A_{1\cup 2}=pa_{G_1}(b_1)\cup pa_{G_2}(b_1)$ and $A_{comp}=(A_{1\cup 2})^c$. Note that: $A_{1\backslash 2}$, $A_{2\backslash 1}\neq \emptyset$, $A_{1\cap 2}$ and $A_{comp}$ can be $\emptyset$; $A_{1\cap 2}$, $A_{1\backslash 2}$, $A_{2\backslash 1}$, and $A_{comp}$ is a partition of $N$.

Consider a function $w$ similar with part (1) that can also be considered as a vector such that the permutations of coordinates in $w$ and in characteristic imsets coincide. 
\begin{itemize}
\item
If $|A_{1\backslash 2}|>1$ and $|A_{2\backslash 1}|>1$, we define $w$ as:
\small{\[
 w(T) = \left\{ \begin{array}{lcl}
c   & \mbox{for} & T=a_i b_j,\ a_i\in A_{1\cap 2} \\
-c & \mbox{for} & T=a_i b_j,\ a_i\notin A_{1\cap 2} \\
-2c & \mbox{for} & T=A_{1\backslash 2}\cup A_{2\backslash 1}\cup\{b_1\} \\
(|A_{1\backslash 2}|+1)\cdot c & \mbox{for} & T=A_{1\backslash 2}\cup\{b_1\} \\
(|A_{2\backslash 1}|+1)\cdot c & \mbox{for} & T=A_{2\backslash 1}\cup\{b_1\} \\
0   & \mbox{for} & \mbox{other} \ T\subset N, |T|>2
\end{array}\right.
\]} \normalsize
where c is a positive number.

Then $\forall \ G_3 \in \mathcal{G}_{m,1}$, we have:
\small{\[
\begin{array}{ccl}
w\cdot c_{G_3} & = & |pa_{G_3}(b_1)\cap A_{1\cap 2}|\cdot c-|pa_{G_3}(b_1)\cap A_{1\backslash 2}|\cdot c\\
 &     &-|pa_{G_3}(b_1)\cap A_{2\backslash 1}|\cdot c
-|pa_{G_3}(b_1)\cap A_{comp}|\cdot c\\
 &     & +(|A_{1\backslash 2}|+1)\cdot c \cdot c_{G_3}(A_{1\backslash 2}\cup\{b_1\})
+(|A_{2\backslash 1}|+1)\cdot c \cdot c_{G_3}(A_{2\backslash 1}\cup\{b_1\}) \\
 &     & -2c\cdot c_{G_3}(A_{1\backslash 2}\cup A_{2\backslash 1}\cup\{b_1\}) \\
 & = & |pa_{G_3}(b_1)\cap A_{1\cap 2}|\cdot c \\
 &     & -|pa_{G_3}(b_1)\cap A_{1\backslash 2}|\cdot c +(|A_{1\backslash 2}|+1)\cdot c \cdot c_{G_3}(A_{1\backslash 2}\cup\{b_1\}) \\
 &     & -|pa_{G_3}(b_1)\cap A_{2\backslash 1}|\cdot c +(|A_{2\backslash 1}|+1)\cdot c \cdot c_{G_3}(A_{2\backslash 1}\cup\{b_1\}) \\
 &     & -2c\cdot c_{G_3}(A_{1\backslash 2}\cup A_{2\backslash 1}\cup\{b_1\}) \\
 &     & -|pa_{G_3}(b_1)\cap A_{comp}|\cdot c
\end{array}
\]} \normalsize
In this equation:
\begin{itemize}
\item
$|pa_{G_3}(b_1)\cap A_{1\cap 2}|\cdot c \leq |A_{1\cap 2}|\cdot c$, where ``=" holds if and only if $A_{1\cap 2}\subset pa_{G_3}(b_1)$;
\item
$-|pa_{G_3}(b_1)\cap A_{1\backslash 2}|\cdot c +(|A_{1\backslash 2}|+1)\cdot c \cdot c_{G_3}(A_{1\backslash 2}\cup\{b_1\}) \leq c$, where ``=" holds if and only if $A_{1\backslash 2}\subset pa_{G_3}(b_1)$;
\item
$-|pa_{G_3}(b_1)\cap A_{2\backslash 1}|\cdot c +(|A_{2\backslash 1}|+1)\cdot c \cdot c_{G_3}(A_{2\backslash 1}\cup\{b_1\}) \leq c$, where ``=" holds if and only if $A_{2\backslash 1}\subset pa_{G_3}(b_1)$;
\item
$-2c\cdot c_{G_3}(A_{1\backslash 2}\cup A_{2\backslash 1}\cup\{b_1\})\leq 0$, where ``=" holds if and only if $(A_{1\backslash 2}\cup A_{2\backslash 1})\nsubseteq pa_{G_3}(b_1)$;
\item
$-|pa_{G_3}(b_1)\cap A_{comp}|\cdot c \leq 0$, where ``=" holds if and only if $pa_{G_3}(b_1)\cap A_{comp} = \emptyset$.
\end{itemize}
The above conditions cannot be satisfied simultaneously, but notice that:
\begin{itemize}
\item
when $pa_{G_3}(b_1)=A_{1\cap 2}$, $w\cdot c_{G_3} =|A_{1\cap 2}|\cdot c+0+0+0+0=|A_{1\cap 2}|\cdot c$;
\item
when $pa_{G_3}(b_1)=A_1$, i.e. $G_3=G_1$, $w\cdot c_{G_3} =|A_{1\cap 2}|\cdot c+c+0+0+0=(|A_{1\cap 2}|+1)\cdot c$;
\item
when $pa_{G_3}(b_1)=A_2$, i.e. $G_3=G_2$, $w\cdot c_{G_3} =|A_{1\cap 2}|\cdot c+0+c+0+0=(|A_{1\cap 2}|+1)\cdot c$;
\item
when $pa_{G_3}(b_1)=A_{1\cup 2}$, $w\cdot c_{G_3} =|A_{1\cap 2}|\cdot c+c+c-2c+0=|A_{1\cap 2}|\cdot c$.
\end{itemize}
Now it is obvious that $w\cdot c_{G_3} \leq (|A_{1\cap 2}|+1)\cdot c$, where ``=" holds if and only if $G_3=G_1$ or $G_2$.

\item
If only one of $|A_{1\backslash 2}|$ and $|A_{2\backslash 1}|$ is 1. Suppose $|A_{1\backslash 2}|=1$ and $|A_{2\backslash 1}|>1$. We define $w$ as:
\small{\[
 w(T) = \left\{ \begin{array}{lcl}
c   & \mbox{for} & T=a_i b_j,\ a_i\in A_1 \\
-c & \mbox{for} & T=a_i b_j,\ a_i\notin A_1 \\
-2c & \mbox{for} & T=A_{1\backslash 2}\cup A_{2\backslash 1}\cup\{b_1\} \\
(|A_{2\backslash 1}|+1)\cdot c & \mbox{for} & T=A_{2\backslash 1}\cup\{b_1\} \\
0   & \mbox{for} & \mbox{other} \ T\subset N, |T|>2
\end{array}\right.
\]} \normalsize
where c is a positive number.

Then $\forall \ G_3 \in \mathcal{G}_{m,1}$, we have:
\small{\[
\begin{array}{ccl}
w\cdot c_{G_3} & = & |pa_{G_3}(b_1)\cap A_{1\cap 2}|\cdot c+|pa_{G_3}(b_1)\cap A_{1\backslash 2}|\cdot c\\
 &     &-|pa_{G_3}(b_1)\cap A_{2\backslash 1}|\cdot c
-|pa_{G_3}(b_1)\cap A_{comp}|\cdot c\\
 &     & +(|A_{2\backslash 1}|+1)\cdot c \cdot c_{G_3}(A_{2\backslash 1}\cup\{b_1\}) 
-2c\cdot c_{G_3}(A_{1\backslash 2}\cup A_{2\backslash 1}\cup\{b_1\}) \\
 & = & |pa_{G_3}(b_1)\cap A_{1\cap 2}|\cdot c \\
 &     & +|pa_{G_3}(b_1)\cap A_{1\backslash 2}|\cdot c \\
 &     & -|pa_{G_3}(b_1)\cap A_{2\backslash 1}|\cdot c +(|A_{2\backslash 1}|+1)\cdot c \cdot c_{G_3}(A_{2\backslash 1}\cup\{b_1\}) \\
 &     & -2c\cdot c_{G_3}(A_{1\backslash 2}\cup A_{2\backslash 1}\cup\{b_1\}) \\
 &     & -|pa_{G_3}(b_1)\cap A_{comp}|\cdot c.
\end{array}
\]} \normalsize
In this equation:
\begin{itemize}
\item
$|pa_{G_3}(b_1)\cap A_{1\cap 2}|\cdot c \leq |A_{1\cap 2}|\cdot c$, where ``='' holds if and only if $A_{1\cap 2}\subset pa_{G_3}(b_1)$;
\item
$|pa_{G_3}(b_1)\cap A_{1\backslash 2}|\cdot c \leq c$, where ``='' holds if and only if $A_{1\backslash 2}\subset pa_{G_3}(b_1)$;
\item
$-|pa_{G_3}(b_1)\cap A_{2\backslash 1}|\cdot c +(|A_{2\backslash 1}|+1)\cdot c \cdot c_{G_3}(A_{2\backslash 1}\cup\{b_1\}) \leq c$, where ``='' holds if and only if $A_{2\backslash 1}\subset pa_{G_3}(b_1)$;
\item
$-2c\cdot c_{G_3}(A_{1\backslash 2}\cup A_{2\backslash 1}\cup\{b_1\})\leq 0$, where ``='' holds if and only if $(A_{1\backslash 2}\cup A_{2\backslash 1})\nsubseteq pa_{G_3}(b_1)$;
\item
$-|pa_{G_3}(b_1)\cap A_{comp}|\cdot c \leq 0$, where ``='' holds if and only if $pa_{G_3}(b_1)\cap A_{comp} = \emptyset$.
\end{itemize}
The above conditions cannot be satisfied simultaneously, but it is similar with the case of ``$|A_{1\backslash 2}|>1$ and $|A_{2\backslash 1}|>1$'' to show that $w\cdot c_{G_3} \leq (|A_{1\cap 2}|+1)\cdot c$, where ``='' holds if and only if $G_3=G_1$ or $G_2$.

\item
If $|A_{1\backslash 2}|=|A_{2\backslash 1}|=1$, we define $w$ as:
\small{\[
 w(T) = \left\{ \begin{array}{lcl}
c   & \mbox{for} & T=a_i b_j,\ a_i\in A_{1\cup 2} \\
-c & \mbox{for} & T=a_i b_j,\ a_i\notin A_{1\cup 2} \\
-2c & \mbox{for} & T=A_{1\backslash 2}\cup A_{2\backslash 1}\cup\{b_1\} \\
0   & \mbox{for} & \mbox{other} \ T\subset N, |T|>2
\end{array}\right.
\]} \normalsize
where c is a positive number.

Then $\forall \ G_3 \in \mathcal{G}_{m,1}$, we have:
\[
\begin{array}{ccl}
w\cdot c_{G_3} & = & |pa_{G_3}(b_1)\cap A_{1\cap 2}|\cdot c \\
 &     & +|pa_{G_3}(b_1)\cap A_{1\backslash 2}|\cdot c
 +|pa_{G_3}(b_1)\cap A_{2\backslash 1}|\cdot c \\
 &     & -2c\cdot c_{G_3}(A_{1\backslash 2}\cup A_{2\backslash 1}\cup\{b_1\}) \\
 &     & -|pa_{G_3}(b_1)\cap A_{comp}|\cdot c .
\end{array}
\]
In this equation:
\begin{itemize}
\item
$|pa_{G_3}(b_1)\cap A_{1\cap 2}|\cdot c \leq |A_{1\cap 2}|\cdot c$, where ``='' holds if and only if $A_{1\cap 2}\subset pa_{G_3}(b_1)$;
\item
$|pa_{G_3}(b_1)\cap A_{1\backslash 2}|\cdot c \leq c$, where ``='' holds if and only if $A_{1\backslash 2}\subset pa_{G_3}(b_1)$;
\item
$|pa_{G_3}(b_1)\cap A_{2\backslash 1}|\cdot c \leq c$, where ``='' holds if and only if $A_{2\backslash 1}\subset pa_{G_3}(b_1)$;
\item
$-2c\cdot c_{G_3}(A_{1\backslash 2}\cup A_{2\backslash 1}\cup\{b_1\})\leq 0$, where ``='' holds if and only if $(A_{1\backslash 2}\cup A_{2\backslash 1})\nsubseteq pa_{G_3}(b_1)$;
\item
$-|pa_{G_3}(b_1)\cap A_{comp}|\cdot c \leq 0$, where ``='' holds if and only if $pa_{G_3}(b_1)\cap A_{comp} = \emptyset$.
\end{itemize}
The above conditions cannot be satisfied simultaneously, but it is similar with the case of ``$|A_{1\backslash 2}|>1$ and $|A_{2\backslash 1}|>1$'' to show that: $w\cdot c_{G_3} \leq (|A_{1\cap 2}|+1)\cdot c$, where ``='' holds if and only if $G_3=G_1$ or $G_2$.
\end{itemize}

\end{itemize}
\end{proof}

\subsection{Proof of Theorem \ref{thm:diagedge}}

\begin{proof} 
We will prove ``if'' and ``only if'' separately.
\begin{itemize}
\item[(1)]
Prove ``if'' part. 

Suppose $G_1$, $G_2\in \mathcal{G}_{m,n}$, and there exists $b_i\in B$ such that $pa_{G_1}(b_i)\neq pa_{G_2}(b_i)$ and $pa_{G_1}(b_j) = pa_{G_2}(b_j)$, $\forall\ b_j\in B$, $b_j\neq b_i$. We need to prove $G_1$ and $G_2$ are neighbors.

Consider an arbitrary graph $G_3\in \mathcal{G}_{m,n}$. We need to prove: $\exists$ a cost vector $w$ such that $w\cdot c_{G_1}=w\cdot c_{G_2}\geq w\cdot c_{G_3}$, where ``='' holds if and only if $G_3=G_1$ or $G_2$.

Define the following graphs (a graphical example will be given in Remark \ref{rem:diagproof}):
\begin{itemize}
\item
$G_1'$, $G_2'$, $G_3' \in \mathcal{G}_{m,1}$ with symptom $B_{m,1}=\{b_i\}$ such that $pa_{G_1'}(b_i)=pa_{G_1}(b_i)$, $pa_{G_2'}(b_i)=pa_{G_2}(b_i)$ and $pa_{G_3'}(b_i)=pa_{G_3}(b_i)$;
\item
$G_0$, $G_3''\in \mathcal{G}_{m,(n-1)}$ with symptoms $B_{m,(n-1)}=B\backslash\{b_i\}$ such that $pa_{G_0}(b_j)=pa_{G_1}(b_j)=pa_{G_2}(b_j)$ and $pa_{G_3''}(b_j)=pa_{G_3}(b_j)$, $\forall\ b_j\in B_{m,(n-1)}$.
\end{itemize}
By Remark \ref{rem:diag1}, with a proper permutation of coordinates, we can write the characteristic imsets of $G_1$, $G_2$ and $G_3$ in the form of:
\[
\begin{array}{cll}
c_{G_1}= & (c_{G_1'}, & c_{G_0}) \\
c_{G_2}= & (c_{G_2'}, & c_{G_0}) \\
c_{G_3}= & (c_{G_3'}, & c_{G_3''})
\end{array}
\]
\begin{itemize}
\item
By Lemma \ref{lem:onesymp}, $G_1'$ and $G_2'$ are neighbors, i.e. $\exists$ a cost vector $w_1$ such that $w_1\cdot c_{G_1'}=w_1\cdot c_{G_2'}\geq w_1\cdot c_{G_3'}$, $\forall\ G_3'\in \mathcal{G}_{m,1}$, where ``='' holds if and only if $G_3'=G_1'$ or $G_2'$.

\item
Since $c_{G_0} \in vert(\mathbf P_{\mathcal{G}_{m,(n-1)}, c})$, $\exists$ a cost vector $w_2$ such that $w_2\cdot c_{G_0} \geq w_2 \cdot c_{G_3''}$, $\forall\ G_3''\in \mathcal{G}_{m,(n-1)}$, where ``='' holds if and only if $G_3''=G_0$.
\end{itemize}
Let $w=(w_1\ w_2)$. We have:
\[
\begin{array}{lllll}
w\cdot c_{G_1} & = & w_1\cdot c_{G_1'}+w_2\cdot c_{G_0}  &    &                          \\
                         & = & w_1\cdot c_{G_2'}+w_2\cdot c_{G_0}  & = & w\cdot c_{G_2} \\
                    & \geq & w_1\cdot c_{G_3'}+w_2\cdot c_{G_3''} & = & w\cdot c_{G_3},
\end{array}
\]
where ``='' holds if and only if i) $G_3'=G_1'$ or $G_2'$, and ii) $G_3''=G_0$, i.e. $G_3=G_1$ or $G_2$.

\item[(2)]
Prove ``only if'' part.

Suppose $G_1$, $G_2\in \mathcal{G}_{m,n}$ are neighbors. i.e. $\exists$ a cost vector $w$ such that $w\cdot c_{G_1} = w\cdot c_{G_2} > w\cdot c_{G} $, $\forall\ G\in\mathcal{G}_{m,n}$, $G\neq G_1$, $G_2$. We are going to prove this part by contradiction.

Suppose $\exists \ b_i$, $b_j\in B$ distinct, $pa_{G_1}(b_i)\neq pa_{G_2}(b_i)$ and $pa_{G_1}(b_j) \neq pa_{G_2}(b_j)$.

Define the following graphs (a graphical example will be given in Remark \ref{rem:diagproof}):
\begin{itemize}
\item
$G_1'$, $G_2' \in \mathcal{G}_{m,1}$ with symptom $B_{m,1}=\{b_i\}$ such that $pa_{G_1'}(b_i)=pa_{G_1}(b_i)$ and $pa_{G_2'}(b_i)=pa_{G_2}(b_i)$;
\item
$G_1''$, $G_2'' \in \mathcal{G}_{m,1}$ with symptom $B_{m,1}=\{b_j\}$ such that $pa_{G_1''}(b_j)=pa_{G_1}(b_j)$ and $pa_{G_2''}(b_j)=pa_{G_2}(b_j)$;
\item
$G_1'''$, $G_2'''\in \mathcal{G}_{m,(n-2)}$ with symptoms $B_{m,(n-2)}=B\backslash\{b_i,b_j\}$ such that $pa_{G_1'''}(b_k)=pa_{G_1}(b_k)$ and $pa_{G_2'''}(b_k)=pa_{G_2}(b_k)$, $\forall\ b_k\in B_{m,(n-2)}$;
\item
$G_3\in \mathcal{G}_{m,n}$ is all the same with $G_1$ but $pa_{G_3}(b_i)=pa_{G_2}(b_i)$;
\item
$G_4\in \mathcal{G}_{m,n}$ is all the same with $G_1$ but $pa_{G_4}(b_j)=pa_{G_2}(b_j)$;
\item
$G_5\in \mathcal{G}_{m,n}$ is all the same with $G_2$ but $pa_{G_5}(b_i)=pa_{G_1}(b_i)$ and $pa_{G_5}(b_j)=pa_{G_1}(b_j)$, notice that $G_5$ might be same with $G_1$.
\end{itemize}
Similarly with part (1), with a proper permutation of coordinates, we can write the characteristic imsets of $G_1$, $G_2$, $G_3$, $G_4$ and $G_5$ in the following form:
\[
\begin{array}{clll}
c_{G_1}= & (c_{G_1'}, & c_{G_1''}, & c_{G_1'''}) \\
c_{G_2}= & (c_{G_2'}, & c_{G_2''}, & c_{G_2'''}) \\
c_{G_3}= & (c_{G_2'}, & c_{G_1''}, & c_{G_1'''}) \\
c_{G_4}= & (c_{G_1'}, & c_{G_2''}, & c_{G_1'''}) \\
c_{G_5}= & (c_{G_1'}, & c_{G_1''}, & c_{G_2'''}) 
\end{array}
\]
With the same permutation of coordinates, $w$ can be written as $w=(w_1\ w_2\ w_3)$. Thus we have:
\begin{itemize}
\item
$G_3\neq G_1$ or $G_2$, which implies:
\[
\begin{array}{rlll}
    & w\cdot c_{G_1} & = & w_1 \cdot c_{G_1'} + w_2 \cdot c_{G_1''} + w_3 \cdot c_{G_1'''} \\
> & w\cdot c_{G_3} & = & w_1 \cdot c_{G_2'} + w_2 \cdot c_{G_1''} + w_3 \cdot c_{G_1'''}  \\
\Longrightarrow & w_1\cdot c_{G_1'} & > & w_1 \cdot c_{G_2'};
\end{array}
\]
\item
$G_4\neq G_1$ or $G_2$, which implies:
\[
\begin{array}{rlll}
    & w\cdot c_{G_1} & = & w_1 \cdot c_{G_1'} + w_2 \cdot c_{G_1''} + w_3 \cdot c_{G_1'''} \\
> & w\cdot c_{G_4} & = & w_1 \cdot c_{G_1'} + w_2 \cdot c_{G_2''} + w_3 \cdot c_{G_1'''}  \\
\Longrightarrow & w_2\cdot c_{G_1''} & > & w_2 \cdot c_{G_2''}.
\end{array}
\]
\end{itemize}
There is a contradiction:
\[
\begin{array}{cccll}
    & w\cdot c_{G_2} & = & w_1 \cdot c_{G_2'} + w_2 \cdot c_{G_2''} + w_3 \cdot c_{G_2'''} & \\
    &                          & < & w_1 \cdot c_{G_1'} + w_2 \cdot c_{G_1''} + w_3 \cdot c_{G_2'''} 
    & = w\cdot c_{G_5} \\
\Longrightarrow & w\cdot c_{G_2} & < & w\cdot c_{G_5}. &
\end{array}
\]
Therefore $G_1$ and $G_2$ cannot be neighbors.

\end{itemize}
\end{proof}

\begin{rem}\label{rem:diagproof}
Two graphical examples will be given for a more intuitive view of the proof of Theorem \ref{thm:diagedge}.
\begin{itemize}
\item
Part (1), the proof of ``if'' statement. In Figure \ref{fig:diagnosis_proof1}, $m=4$, $n=3$ and $b_i=b_1$.
\begin{figure}[!htp]
\begin{center}
\scalebox{0.8}{
\includegraphics{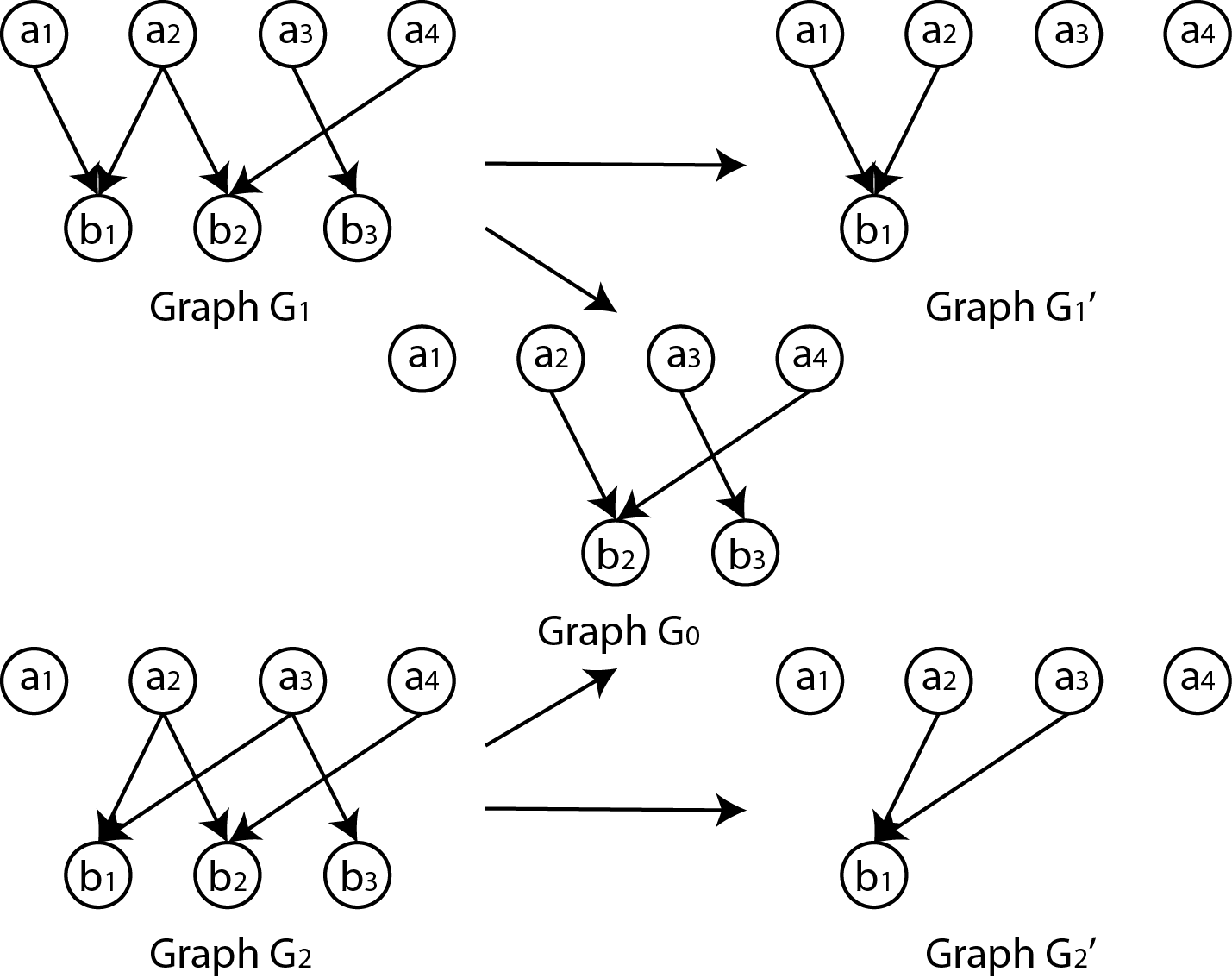}
}
\end{center}
\caption{An example for the proof of Theorem \ref{thm:diagedge}, part (1)}
\label{fig:diagnosis_proof1}
\end{figure}

\item
Part (2), the proof of ``only if'' statement. In Figure \ref{fig:diagnosis_proof2}, $m=4$, $n=3$, $b_i=b_1$ and $b_j=b_2$.
\begin{figure}[!htp]
\begin{center}
\scalebox{0.8}{
\includegraphics{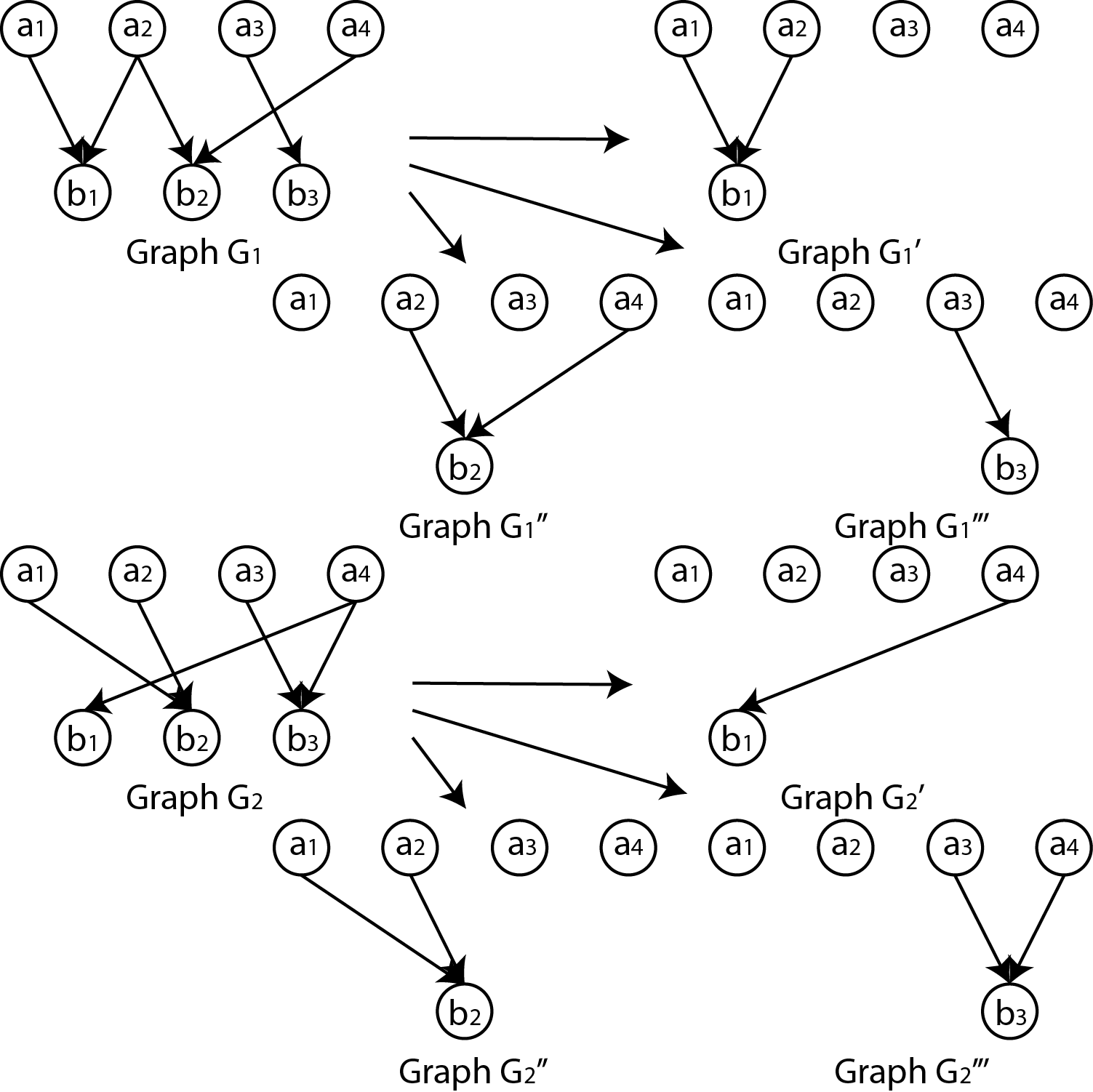}
}
\end{center}
\caption{An example for the proof of Theorem \ref{thm:diagedge}, part (2)}
\label{fig:diagnosis_proof2}
\end{figure}
\end{itemize}
\end{rem}

\subsection{Proof of Theorem \ref{thm:diagnedge}}

\begin{proof} 
By Theorem \ref{thm:diagedge}, $\forall H\in \mathcal{G}_{m,n}$, $G$ and $H$ are neighbors if and only if: $\exists\ b_k\in B$ such that $pa_{G}(b_k)\neq pa_{H}(b_k)$ and $pa_{G}(b_j) = pa_{H}(b_j)$, $\forall\ b_j\in B$ and $b_j\neq b_k$.

Now fix $b_i \in B$. Define graphs:
\begin{itemize}
\item
$G'$, $H' \in \mathcal{G}_{m,1}$ with symptom $B_{m,1}=\{b_i\}$ such that $pa_{G'}(b_i)=pa_{G}(b_i)$ and $pa_{H'}(b_i)=pa_{H}(b_i)$;
\item
$G''$, $H''\in \mathcal{G}_{m,(n-1)}$ with symptoms $B_{m,(n-1)}=B\backslash\{b_i\}$ such that $pa_{G''}(b_j)=pa_{G}(b_j)$ and $pa_{H''}(b_j)=pa_{H}(b_j)$, $\forall\ b_j\in B_{m,(n-1)}$.
\end{itemize}
Since $G$ and $H$ are neighbors and $G'\neq H'$ will lead to $G'' = H''$, and by Proposition \ref{prop:diag3} there are $2^m$ graphs in $\mathcal{G}_{m,1}$, there are $2^m-1$ different choices of $H'$s, and each corresponds to a different neighbor of $G$.

We can use the same strategy for every $b_i \in B$, i.e. we can find $2^m-1$ neighbors from each fixed $b_i\in B$. It is easy to see that these neighbors are all distinct:
if $H_1$, $H_2$ are all the same with $G$ but $pa_{G}(b_i)\neq pa_{H_1}(b_i)$ and  $pa_{G}(b_j)\neq pa_{H_2}(b_j)$, where $b_i$, $b_j\in B$ are distinct, then this implies $pa_{H_2}(b_i) = pa_{G}(b_i)\neq pa_{H_1}(b_i)$, i.e. $H_1$ and $H_2$ are different.
Therefore the total number of neighbors for $G$ is: $n\cdot (2^m-1)$.
\end{proof}

\section{Prove $\mathbf P_{m,n}$ is simple using linear algebra}\label{proof:simple}

Recall that in Section \ref{BN:diagnosis}, we first proved that $\mathbf P_{m,1}$ is a simplex $\Delta_{2^m-1}$, and then we proved that $\mathbf P_{m,n}$ is a direct product of $n$ many $\Delta_{2^m-1}$, which implies that $\mathbf P_{m,n}$ is a simple polytope with dimension $n \cdot (2^m-1)$. Now we are going to show another flow to prove that $\mathbf P_{m,n}$ is simple.

First, we will use linear algebra to show that $\mathbf P_{m,n}$ has dimension $n \cdot (2^m-1)$.
We adopt the notation from Section \ref{BN:diagnosis}. Given $N$, by Proposition \ref{prop:diag1} and Proposition \ref{prop:diag4}, we can define $\mathcal{S}_{m,n}$ as the support of $\{c_G: G\in\mathcal{G}_{m,n}\}$, i.e.:
\[
\mathcal{S}_{m,n}=\{T:\exists\ G\in\mathcal{G}_{m,n} \mbox{ such that } c_G(T)=1\}
\subset \mathcal{P}(N),
\]
where $\mathcal{P}(N)$ is the power set of $N$.

\begin{thm}\label{thm:diagdim}
Fix $m$ and $n$. The dimension of $\mathbf P_{m,n}$ is exactly $n\cdot (2^m-1)$.
\end{thm}

\begin{proof}
Similar with imsets, we can consider the standard basis $\mathbf{e_T}$, $T\subset N$, as functions $\mathbf{e_T} : \mathcal P(N) \mapsto \mathbb Z$ such that $\forall \ T_0 \subset N$, $\mathbf{e_T}(T_0)=1$ if $T_0=T$, and $0$ otherwise. Each $\mathbf{e_T}$ can also be considered as a vector with coordinates $T_0 \subset N$.

It is obvious that: 1 ) $\{c_G$, $G\in \mathcal{G}_{m,n}\}\subset \mathbb{R}^{2^{m+n}-(m+n+1)}$; 2 ) $\{\mathbf{e_T}, T\in \mathcal{S}_{m,n}\}$ is a basis of $\mathbb{R}^{n\cdot (2^m-1)}$ that is embedded in $\mathbb{R}^{2^{m+n}-(m+n+1)}$ (Proposition \ref{prop:diag4}); and 3 )$\{c_G$, $G\in \mathcal{G}_{m,n}\}$ can be written as a linear combination of $\{\mathbf{e_T}, T\in \mathcal{S}_{m,n}\}$. We are going to prove that $\{\mathbf{e_T}, T\in \mathcal{S}_{m,n}\}$ can be expressed as a linear combination of $\{c_G$, $G\in \mathcal{G}_{m,n}\}$. Notice that $\{\mathbf{e_T}, T\in \mathcal{S}_{m,n}\}$ is equivalent with $\{\mathbf{e_T}, T\subset N$ and $T$ has the form of $a_{i_1} \ldots a_{i_k} b_j$, where $1\leq k \leq m$, $\{i_1,\ldots,i_k\}\subseteq \{1,\ldots,m\}$ and
$j\in \{1,\ldots,n\}\}$ (Proposition \ref{prop:diag1}), we can prove the statement by induction on $|T|$.

\begin{itemize}
\item When $|T|=2$ (i.e. $k=1$), i.e. $T=a_i b_j$, where $a_i\in A$ and $b_j\in B$, we know $c_G=\mathbf{e_T}$, where $G \in \mathcal G_{m,n}$ has only one edge $a_i \rightarrow b_j$.

\item Suppose $\forall\ T$, $T$ has the form in Proposition \ref{prop:diag1} and $|T| \leq k$, $\mathbf{e_T}$ can be written as a linear combination of $\{c_G$, $G\in \mathcal{G}_{m,n}\}$. Now consider $T_k=a_{i_1} \ldots a_{i_k} b_j$, where $\{i_1,\ldots,i_k\}\subseteq \{1,\ldots,m\}$ and $j\in \{1,\ldots,n\}$. \\
Let $G \in \mathcal G_{m,n}$ have $k$ edges: $a_{i_l} \rightarrow b_j$, $l=1\ldots k$. Then:
\[
\mathbf{e_{T_k}}=c_G-\sum\limits_{T_a\subset \{a_{i_1}, \ldots, a_{i_k}\},0<|T_a|<k} \mathbf{e_{T_a\cup \{b_j\}}}.
\]
Since $\forall\ T_a\subset \{a_{i_1}, \ldots, a_{i_k}\}$, $0<|T_a|<k$ (i.e. $T_a\subsetneq \{a_{i_1}, \ldots, a_{i_k}\}$), $|T_a\cup b_j| \leq k$, $\mathbf{e_{T_a\cup b_j}}$ can be expressed as a linear combination of $\{c_G$, $G\in \mathcal{G}_{m,n}\}$. Therefore, $\mathbf{e_{T_k}}$ can be written as a linear combination of $\{c_G$, $G\in \mathcal{G}_{m,n}\}$.
\end{itemize}

\end{proof}

A special case of $n=1$ in Theorem \ref{thm:diagdim} and Proposition \ref{prop:diag3} claims that $\mathbf{P}_{m,1}$ has $2^m$ vertices and dimension $2^m-1$. This directly lead to Corollary \ref{cor:diagsimplex}.

\begin{cor}\label{cor:diagsimplex}
Fix m, $\mathbf{P}_{m,1}$ is a simplex with dimension $2^m-1$, i.e. $\mathbf{P}_{m,1}=\Delta_{2^m-1}$.
\end{cor}

Lemma \ref{lem:onesymp} is an immediate result of Corollary \ref{cor:diagsimplex}, while Theorem \ref{thm:diagnedge} and Theorem \ref{thm:diagprod} can be obtained based on Lemma \ref{lem:onesymp} and Corollary \ref{cor:diagsimplex} using the same proofs in Section \ref{BN:diagnosis}. It is worth mentioning that Theorem \ref{thm:diagdim} and Theorem \ref{thm:diagnedge} imply that $\mathbf{P}_{m,n}$ is a simple polytope with dimension $n \cdot (2^m-1)$ because the number of neighbors for each vertex equals to the dimension of the polytope. In 2000, V. Kaibel and M. Wolff proved that a zero-one polytope is simple if and only if it equals to a direct product of zero-one simplices \cite{Kaibel}. Recall that cim-polytopes are zero-one polytopes, we are able to conclude that $\mathbf{P}_{m,n}$ is a direct product of zero-one simplices \cite{Kaibel}. Our progress is that we proved a even strong result in Theorem \ref{thm:diagprod} with an intuitive graphical interpretation of each simplex in the direct product.

\section{Proofs in Section \ref{BN:generalize}}\label{proof:generalize}

\subsection{Proof of Theorem \ref{thm:genprod}}

\begin{proof} We are going to prove the equality by induction on $n$. Since $n \geq 2$, we start the induction from $n=2$.
\begin{itemize}
\item
$n=2$. It is obvious since there are only two vertices in $\mathbf{P}_{[n]}$: $(1)$ and $(0)$. So $\mathbf{P}_{[n]}$ is a line segment which is a simplex of dimension $1$, i.e. $\mathbf{P}_{[n]} = \Delta_{1}$.
\item
Fix $q\in \mathbb{Z}_+$. Suppose the equality holds for $\mathbf{P}_{[n]}$, $\forall \ n<q$, and we need to prove that it also holds for $\mathbf{P}_{[q]}$. Define notation $N_{[k]}=\{a_{[1]},\ldots, a_{[k]}\}$ for $k=1,\ldots,q$.

First, we want to prove: $\mathbf{P}_{[q]}\subseteq \mathbf{P}_{[q-1]}\times \Delta_{2^{q-1}-1}$.

$\forall\ G\in \mathcal{G}_{[q]}$, we can define graphs:
\begin{itemize}
\item
$G'$ is the induced subgraph of $G$ for $N_{[q-1]}$, which implies $c_{G'} \in \mathbf{P}_{[q-1]}$;
\item
$G''$ is a graph over $N$ such that the only edges in $G''$ are $a_{[i]} \rightarrow a_{[q]}$, where $a_{[i]} \in pa_G(a_{[q]})$. Consider a diagnosis model where $N_{[q-1]}$ is the set of diseases and $a_{[q]}$ is the symptom, then we can see that $c_{G''} \in \mathbf{P}_{q-1,1}=\Delta_{2^{q-1}-1}$.
\end{itemize}
Now, with a proper permutation of coordinates (see Remark \ref{rem:gen1}), we can write $c_G$ in the form of:
\[
c_G=(c_{G'}\ \ c_{G''}).
\]
Since $vert(\mathbf{P}_{[q]}) = \{c_G: G\in \mathcal{G}_{[q]}\}$, $\forall\ x\in\mathbf{P}_{[q]}$, with the same permutation of coordinates, we have:
\begin{equation}\label{eq:genprod1}
x=\sum\limits_{G\in\mathcal{G}_{[q]}} \alpha_G c_G =
 (\sum\limits_{G\in\mathcal{G}_{[q]}} \alpha_G c_{G'} \ ,\  \sum\limits_{G\in\mathcal{G}_{[q]}} \alpha_G c_{G''}),
\end{equation}
where $0 \leq \alpha_G\leq 1$, $\forall\ G\in \mathcal{G}_{[q]}$ and $\sum\limits_{G\in\mathcal{G}_{[q]}} \alpha_G =1$.

Notice that $\sum\limits_{G\in\mathcal{G}_{[q]}} \alpha_G c_{G'} \in \mathbf{P}_{[q-1]}$ and $\sum\limits_{G\in\mathcal{G}_{[q]}} \alpha_G c_{G''}\in \Delta_{2^{q-1}-1}$, Equation \eqref{eq:genprod1} implies $x\in \mathbf{P}_{[q-1]} \times \Delta_{2^{q-1}-1}$.
Hence:
\[
\mathbf{P}_{[q]}\subseteq  \mathbf{P}_{[q-1]} \times \Delta_{2^{q-1}-1}.
\]

Second, we want to prove: $ \mathbf{P}_{[q-1]} \times \Delta_{2^{q-1}-1} \subseteq \mathbf{P}_{[q]}$.

Let $\mathcal{G}_{[q-1]}$ has nodes $N_{[q-1]}$, and $\mathcal{G}_{q-1,1}$ has diseases $N_{[q-1]}$ and symptom $a_{[q]}$. $\forall \ G'\in \mathcal{G}_{[q-1]}$ and $G'' \in \mathcal{G}_{q-1,1}$, we can define $G\in \mathcal{G}_{[q]}$ by extending $G'$ as following: add a node $a_{[q]}$ and edges $(a_{[i]}, a_{[q]})$, $\forall a_{[i]}\in pa_{G''}(a_{[q]})$, to $G'$. 
We can write $c_G$ in the form of $c_G=(c_{G'}\ \ c_{G''})$.

$\forall x\in \mathbf{P}_{[q-1]} \times \Delta_{2^{q-1}-1} $, $x$ can be written as:
\[
\begin{array}{ccl}
x & = &
 (\sum\limits_{G'\in\mathcal{G}_{[q-1]}} \beta_{G'} c_{G'} \ ,\  \sum\limits_{G''\in\mathcal{G}_{q-1,1}} \gamma_{G''} c_{G''})
=\sum\limits_{G'\in\mathcal{G}_{[q-1]}} \sum\limits_{G''\in\mathcal{G}_{q-1,1}}
   \beta_{G'} \gamma_{G''} (c_{G'}\ ,\  c_{G''}) \\
   & = & \sum\limits_{G'\in\mathcal{G}_{[q-1]}} \sum\limits_{G''\in\mathcal{G}_{q-1,1}}
  ( \beta_{G'} \gamma_{G''}) \  c_G\ \ ,
\end{array}
\]
where $0 \leq \beta_{G'}$, $\gamma_{G''}\leq 1$, $\forall G'\in \mathcal{G}_{[q-1]}$, $\forall G''\in \mathcal{G}_{q-1,1}$, and $\sum_{G'\in\mathcal{G}_{[q-1]}} \beta_{G'} =1$, $\sum_{G''\in\mathcal{G}_{q-1,1}} \gamma_{G''} =1$. 

Notice that
\[
\sum\limits_{G'\in\mathcal{G}_{[q-1]}} \sum\limits_{G''\in\mathcal{G}_{q-1,1}}
  ( \beta_{G'} \gamma_{G''}) = \sum\limits_{G'\in\mathcal{G}_{[q-1]}} \beta_{G'}( \sum\limits_{G''\in\mathcal{G}_{q-1,1}}
   \gamma_{G''})= \sum\limits_{G'\in\mathcal{G}_{[q-1]}}\beta_{G'}=1 .
\]
This leads to $x\in \mathbf{P}_{[q]}$. Hence:
\[
\mathbf{P}_{[q-1]} \times \Delta_{2^{q-1}-1} = \mathbf{P}_{[q-1]}\times \mathbf{P}_{q,1} \subseteq \mathbf{P}_{[q]}.
\]
By induction on $n$, we finish the proof by:
\begin{center}
$\mathbf{P}_{[q]} = \mathbf{P}_{[q-1]}\times \mathbf{P}_{q-1,1}= (\Delta_{2^1-1} \times  \cdots \times \Delta_{2^{q-2}-1} ) \times \Delta_{2^{q-1}-1} =  \Delta_{2^1-1} \times  \cdots \times \Delta_{2^{q-1}-1} $.
\end{center}

\end{itemize}
\end{proof}

\subsection{Proof of Theorem \ref{thm:genedge}}

\begin{proof} 
The proof from the view of graph theory will be very similar with the proof of Theorem \ref{thm:diagedge}, so we are going to give a proof from the view of polyhedral geometry, i.e. prove that: ``$\exists$ vertices of $v^1$, $v^2 \in \mathbf{P}_{[n]}$ such that $\mathbf x = \beta v^1 + (1-\beta ) v^2$ where $0 \leq \beta \leq 1$, and $v^1$, $v^2$ form an edge in $\mathbf P_{[n]}$'' if and only if ``$\mathbf x$ can be written in the form of $\mathbf x=(v_1,\ldots,v_{i-1},e_i,v_{i+1},\ldots,v_{n-1})$, $i\in\{1,\ldots, n-1\}$''.

We will prove ``if'' and ``only if'' separately.
\begin{itemize}
\item[(1)]
Prove ``if'' part. 

Suppose $\mathbf x$ has the form $\mathbf x=(v_1,\ldots,v_{i-1},e_i,v_{i+1},\ldots,v_{n-1})$.

Since $e_i$ belongs to an edge on $\Delta_{2^i-1}$, we can find two vertices $v_i^1$, $v_i^2 \in \Delta_{2^i-1}$ which form this edge, and this implies $e_i = \beta v_i^1 + (1-\beta v_i^2)$, $0 \leq \beta \leq 1$. Suppose the cost vector for this edge is $w_i^e$, then for any $v_i^3  \in vert(\Delta_{2^i-1})$, $w_i^e v_i^3 \leq w_i^e v_i^1 = w_i^e v_i^2$, where ``='' holds if and only if $v_i^3 = v_i^1$ or $v_i^3 = v_i^2$.

We can also find $w_j^v$ which is a cost vector for vertex $v_j$ in $\Delta_{2^j-1}$, $j\in \{1, \ldots, n-1\}\backslash \{i\}$. Still, we have: $\forall v_j^3  \in vert(\Delta_{2^j-1})$, $w_j^v v_j^3 \leq w_j^v v_j$, where ``='' holds if and only if $v_j^3 = v_j$.

Now let $v^1 = (v_1, \ldots, v_{i-1}, v_i^1, v_{i+1}, \ldots, v_{n-1})$, $v^2 = (v_1, \ldots, v_{i-1}, v_i^2, v_{i+1}, \ldots, v_{n-1})$ and $w = (w_1^v, \ldots, w_{i-1}^v, w_i^e, w_{i+1}^v, \ldots, w_{n-1}^v)$. Obviously $\mathbf x = \beta v^1 + (1-\beta ) v^2$, where $0 \leq \beta \leq 1$. In addition, $\forall v^3=(v_1^3, \ldots, v_{n-1}^3) \in vert(\mathbf P_{[n]})$, we have:
\[
\begin{array}{lclcl}
w v^3 = w_i^e v_i^3 + \sum\limits_{j=1,\ j\neq i}^{n-1} w_j^v v_j^3  & \leq & w_i^e v_i^1 + \sum\limits_{j=1,\ j\neq i}^{n-1} w_j^v v_j  &  =  &  w v^1    \\
                                                                                                               &   =   & w_i^e v_i^2 + \sum\limits_{j=1,\ j\neq i}^{n-1} w_j^v v_j  &  =  &  w v^2,
\end{array}
\]
where ``='' holds if and only if $v^3=v^1$ or $v^3=v^2$, i.e. $v^1$ and $v^2$ form an edge on $\mathbf P_{[n]}$.

\item[(2)]
Prove ``only if'' part.

Suppose $\exists \ v^1=(v_1^1, \ldots v_{n-1}^1)$, $v^2=(v_1^2, \ldots v_{n-1}^2) \in vert(\mathbf{P}_{[n]})$ such that $\mathbf x = \beta v^1 + (1-\beta ) v^2$ where $0 \leq \beta \leq 1$, and $v^1$, $v^2$ form an edge in $\mathbf P_{[n]}$. If we can prove that $\exists \ i\in \{1, \ldots, n-1\}$ such that $v_i^1\neq v_i^2$ and $v_j^1=v_j^2$, $\forall \ j\in \{1, \ldots, n-1\}\backslash \{i\}$, then $\mathbf x$ has the form $\mathbf x=(v_1,\ldots,v_{i-1},e_i,v_{i+1},\ldots,v_{n-1})$, where $e_i$ is on the edge of $\Delta_{2^i-1}$ formed by $v_i^1$ and $v_i^2$. We are going to prove this statement by contradiction.

Suppose $\exists \ i$, $j \in \{1, \ldots, n-1\}$ distinct such that $v_i^1 \neq v_i^2$ and $v_j^1 \neq v_j^2$, but $v^1$ and $v^2$ still form an edge on $\mathbf P_{[n]}$. Let $w=(w_1, \ldots, w_{n-1})$ be the cost vector for this edge, i.e. $\forall v^3 =(v_1^3, \ldots v_{n-1}^3) \in vert(\mathbf{P}_{[n]})$, $w v^3 \leq w v^1 = w v^2$ where ``='' holds if and only if $v^3=v^1$ or $v^3=v^2$.
\begin{itemize}
\item
If we set $v^3$ as following: $v_i^3 = v_i^2$, $v_k^3 = v_k^1$, $\forall \ k\in \{1,\ldots, n-1\} \backslash \{i\}$. Obviously $v^3 \neq v^1$ and $v^3 \neq v^2$. Thus:
\[
\begin{array}{rcl}
w v^3 = w_i v_i^2 + \sum\limits_{k=1,\ k\neq i}^{n-1} w_k v_k^1  & < & w v^1 = \sum\limits_{k=1}^{n-1} w_k v_k^1 =  w_i v_i^1 + \sum\limits_{k=1, \ k\neq i}^{n-1} w_k v_k^1    \\
                                                 \Longrightarrow \ \ \ \ \  w_i v_i^2  & < & w_i v_i^1.
\end{array}
\]

\item
If we set $v^3$ as following: $v_j^3 = v_j^2$, $v_k^3 = v_k^1$, $\forall \ k\in \{1,\ldots, n-1\} \backslash \{j\}$. Obviously $v^3 \neq v^1$ and $v^3 \neq v^2$. Thus:
\[
\begin{array}{rcl}
w v^3 = w_j v_j^2 + \sum\limits_{k=1,\ k\neq j}^{n-1} w_k v_k^1  & < & w v^1 = \sum\limits_{k=1}^{n-1} w_k v_k^1 =  w_j v_j^1 + \sum\limits_{k=1, \ k\neq j}^{n-1} w_k v_k^1    \\
                                                 \Longrightarrow \ \ \ \ \  w_j v_j^2  & < & w_j v_j^1.
\end{array}
\]
\end{itemize}

Now we set $v^3$ as following: $v_i^3 = v_i^1$, $v_j^3 = v_j^1$, $v_k^3 = v_k^2$, $\forall \ k\in \{1,\ldots, n-1\} \backslash \{i, j\}$. Then we have:
\[
w v^3 = w_i v_i^1 + w_j v_j^1 + \sum\limits_{k=1,\ k\neq i, j}^{n-1} w_k v_k^2  >  
w_i v_i^2 + w_j v_j^2 + \sum\limits_{k=1,\ k\neq i, j}^{n-1} w_k v_k^2  =
\sum\limits_{k=1}^{n-1} w_k v_k^2 = w v^2,
\]
i.e. $w v^3 > w v^2$, which is a contradiction with our assumption. 

\end{itemize}
\end{proof}


\section{Discussion}\label{dis}

\subsection{Connection to K2 algorithm}
If we consider a criterion $\mathcal Q$ which is a regular criterion, then it was proved that $\mathcal Q$ can be written as $\mathcal Q(G,D)=s(D)-\langle r_D, c_G\rangle$, where the entropy $s(D)$ and the data vector $r_D$ only depends on the data $D$ and $c_G$ is the characteristic imset of $G$ \cite{studeny2008, lindner2012}. 
Once the cim-polytope can be written as a direct product of a sequences of simplices, we are able to find the optimal BN structure by maximizing a target function in each simplex: given data
$D \in DATA(N,d)$, 
\begin{align}\label{eq:k2}
\max_{G \in \mathcal G_{[n], \Omega}} \mathcal Q(G,D)
\Longrightarrow
\min_{\mathbf x \in \mathbf P_{\mathcal G_{[n], \Omega}, c}} r_D^T \mathbf x
= \sum\limits_{i=2}^n \min_{\mathbf {x_i} \in \Delta_{2^{|\Omega_{i}|}-1}} r_{D,i}^T \mathbf {x_i} ,
\end{align}
where $\mathbf {x_i}$ contains the coordinates $\{T \subseteq \Omega_i \cup \{a_{[i]}\}: |T| \geq 2$, $ a_{[i]} \in T$, $ a_{[j]} \notin T$, $\forall j>i\}$ in $\mathbf x$, and the coordinates of $r_{D,i}^T$ matches the coordinates of $\mathbf {x_i}$.
This implies that we can find the optimal parent sets of $a_{[i]}$, $i = 2, \ldots, n$, sequentially until we obtain the whole BN structure, which will be exactly the optimal BN structure in $\mathcal G_{[n], \Omega}$.

Equation \eqref{eq:k2} gives a polyhedral geometric insight of the K2 algorithm \cite{cooper1992}, which is a well-known heuristic method in learning Bayesian networks. Recall that in K2 algorithm, an ordering on the nodes is also fixed and parent sets of $a_{[i]}$, $i = 2, \ldots, n$, are also determined sequentially. However, in order to find the optimal BN, Equation \eqref{eq:k2} claims that we need to find $G_i \in \mathcal G_{|\Omega_i|,1}$ such that $r_{D,i}^T c_{G_i} = \min_{\mathbf {x_i} \in \Delta_{2^{|\Omega_{i}|}-1}} r_{D,i}^T \mathbf {x_i}$, while the K2 algorithm obtain each parent set $pa_G(a_{[i]})$ by adding nodes to $\emptyset$ stepwisely (or removing nodes from $\{a_{[1]}, \ldots, a_{[i-1]}\}$ stepwisely), which cannot guarantee that the resulting parent sets are optimal (see Example \ref{ex:k2counter} for a counter-example).

\begin{ex}\label{ex:k2counter} 
Consider $\mathcal G_{3,1}$. The characteristic imsets of all possible graphs in $\mathcal G_{3,1}$ is listed as a matrix:
\footnotesize{
$$ \bordermatrix{ & \cr
	&  c_{G_0} \cr	&  c_{G_1} \cr	&  c_{G_2} \cr
	&  c_{G_3} \cr	&  c_{G_{12}} \cr	&  c_{G_{23}} \cr
	&  c_{G_{13}} \cr	&  c_{G_{123}} } 
= \bordermatrix{ T & a_1 b_1 & a_2 b_1 & a_3 b_1 & a_1 a_2 b_1 & a_1 a_3 b_1 & a_2 a_3 b_1 & a_1 a_2 a_3 b_1 \cr
	       &  0  &  0  &  0 & 0 & 0 & 0 & 0  \cr
	       &  1  &  0  &  0 & 0 & 0 & 0 & 0  \cr
              &  0  &  1  &  0 & 0 & 0 & 0 & 0  \cr
	       &  0  &  0  &  1 & 0 & 0 & 0 & 0  \cr
	       &  1  &  1  &  0 & 1 & 0 & 0 & 0  \cr
	       &  1  &  0  &  1 & 0 & 1 & 0 & 0  \cr
	       &  0  &  1  &  1 & 0 & 0 & 1 & 0  \cr
	       &  1  &  1  &  1 & 1 & 1 & 1 & 1  
}$$
}\normalsize
We are going to give counter-examples that the resulting BN of the K2 algorithm is not the optimal solution.
\begin{itemize}
\item
Forward selection, i.e. each parent set $pa_G(_{[i]})$ is obtained by adding nodes to $\emptyset$ stepwisely. Suppose $r_D^T=(-1, -2, -1, -3, -10, -4, 20)$ which satisfies $r_D^T c_{G_{13}} = -12 < r_D^T c_G$, $\forall G\in \mathcal G_{3,1}$, $G \neq G_{13}$, i.e. the optimal graph is $G_{13}$. 
In K2 algorithm, we start from $pa_G(b_1)=\emptyset$. Next, $a_2$ is added to $pa_G(b_1)$ because $r_D^T c_{G_{2}}=-2 < r_D^T c_{G_{1}}= r_D^T c_{G_{3}}= -1$. Then $a_3$ is added to $pa_G(b_1)$ because $r_D^T c_{G_{23}}=-7 < r_D^T c_{G_{12}}= -6$. Procedure ends here because $r_D^T c_{G_{23}}=-7 < r_D^T c_{G_{123}}= -1$. The graph chosen by K2 algorithm, $G_{23}$, is not the optimal graph.

\item
Backward selection, i.e. each parent set $pa_G(a_{[i]})$ is obtained by removing nodes from $\{a_{[1]}, \ldots, a_{[i-1]}\}$ stepwisely. Suppose $r_D^T=(-3, -1, -1, 3, 3, 0, 10)$ which satisfies $r_D^T c_{G_{1}} = -3 < r_D^T c_G$, $\forall G\in \mathcal G_{3,1}$, $G \neq G_{1}$, i.e. the optimal graph is $G_{1}$. 
In K2 algorithm, we start from $pa_G(b_1)=\{a_1,a_2,a_3\}$. Next, $a_1$ is removed from $pa_G(b_1)$ because $r_D^T c_{G_{23}}=-2 < r_D^T c_{G_{12}}= r_D^T c_{G_{13}}= -1$. Procedure ends here because $r_D^T c_{G_{23}}=-2 < r_D^T c_{G_{2}}= r_D^T c_{G_{3}} =-1$. The graph chosen by K2 algorithm, $G_{23}$, is not the optimal graph.
\end{itemize}
\end{ex}


\subsection{Open problems}\label{diss:BNopen}

Further work and open problems are still left in this topic.
As we mentioned before, the main purpose of studying the structure of cim-polytopes is reducing the time complexity of learning Bayesian networks by suggesting polyhedral geometry techniques. But the reality is that even we have simplified our problem of learning BNs to LP problems over each simplex (see Equation \eqref{eq:k2}) in the direct produce showed in Theorem \ref{thm:genprod} and Equation \eqref{eq:genprod2}, and have described all edges and facets of these simplices (see Section \ref{BN:diagnosis}), if the number of nodes is large, the procedure of searching the optimal solutions in each simplex may still be very time-consuming. In this sense, simulations and analysis on real datasets are necessary to compare the solution and time complexity of our method with other existing classifiers \cite{vomlel2012}. On the other hand, we also need to study on the misspecification (i.e. the underlying ordering of nodes is misspecified)
 and data sensitivity problems of our method via simulations.


Another way to reduce the time complexity is considering setting up a maximum number of parents to control the model complexity, especially when the number of nodes is too large. In this case, since the underlying ordering is fixed, the cim-polytope is still a direct product of simplices. Thus all edges of the cim-polytope can be found similarly with Theorem \ref{thm:genedge}, but the expression of facets for each simplex is not clear.

Notice that all conclusion and discussion in this paper until now are all based on a fixed underlying ordering of nodes. However, in practice, it is often hard to decide such an ordering. One way to compromise is that we can fix the ordering of some of the nodes, and consider every permutation of the rest nodes. For instance, when we use SNP data to examine phenotypes, we are more interested in how genes affect phenotypes and how phenotypes affect each other. Thus we can consider DAGs where all edges between SNPs and edges from phenotypes to SNPs are forbidden, i.e. we only need to consider the permutation of phenotypes.


This paper focuses on the case that all random variables in $N$ are finite random variables. It is still an open problem that how to generalize our method to the case that some or all of the random variables in $N$ are continuous random variables.

\section*{Acknowledgment}

The authors would like to thank Drs.~R. Hemmecke, M. Studen\'y, and B. Sturmfels for
their useful advise. 

\bibliographystyle{siam}
\bibliography{imset}

\end{document}